\documentclass[12pt]{amsart}

\usepackage{amsmath,amssymb,euscript,fullpage,tikz,hyperref,eufrak,verbatim}
\usetikzlibrary{decorations.markings,decorations.pathreplacing,calc,arrows}
\tikzset{
doublearrow/.style={decoration={
  markings,
  mark=at position .6 with {\arrow{>}}},postaction={decorate}, draw, thick, double distance=2pt,>= stealth}
}
\usepackage[hmargin=1 in, vmargin = 0.5 in]{geometry}
\usepackage[all]{xy}

\allowdisplaybreaks[1]

\newtheorem{theorem}{Theorem}
\newtheorem{lemma}[theorem]{Lemma}
\newtheorem{proposition}[theorem]{Proposition}
\newtheorem{corollary}[theorem]{Corollary}
\newtheorem{conjecture}[theorem]{Conjecture}

\newtheorem{definition}[theorem]{Definition}
\newtheorem{proofs}[theorem]{Proof}

\theoremstyle{definition}
\newtheorem{remark}{Remark}
\numberwithin{theorem}{subsection}

\newcommand{\aph}{\alpha}

\newcommand{\spn}{\mathfrak{sp}}

\newcommand{\ejj}[1][j]{[e_{1}[e_1[\ldots[e_{#1-1}e_{#1}]\ldots]}
\newcommand{\eij}[3]{[e_{#1}[\ldots[e_{#2}e_{#3}]\ldots]}
\newcommand{\eiaj}[2]{[e_{#1}[\ldots[e_1[e_1\ldots[e_{#2-1}e_{#2}]\ldots] }
\newcommand{\eia}[2]{[e_{#1}[\ldots[e_1[e_1\ldots[e_{#2}e_{#2+1}]\ldots] }
\newcommand{\eiai}[1]{[e_{#1+1}[\ldots[e_1[e_1\ldots[e_{#1+1}e_{#1+2}]\ldots] }
\newcommand{\bij}[2]{[e_{#1}e_{#2}]}

\newcommand{\beij}[3]{[\bar{e}_{#1}[\ldots,[\bar{e}_{#2}\bar{e}_{#3}]\ldots]}
\newcommand{\bbij}[2]{[\bar{e}_{#1},\bar{e}_{#2}]}
\newcommand{\be}{\bar{e}}

\newcommand{\ee}{[e_1[e_1[e_2e_3]]]}
\newcommand{\eijk}[3]{[e_{#1}[e_{#2}e_{#3}]]}
\newcommand{\eijkr}[4]{[e_{#1}[e_{#2}[e_{#3}e_{#4}]]]}

\newcommand{\eiajj}[4]{[e_{#1}[e_{#2}[\ldots[e_1[e_1\ldots[e_{#3}e_{#4}]\ldots] }
\newcommand{\eiajfive}[5]{[e_{#1}[e_{#2}[\ldots[e_1[e_1\ldots[e_{#3}[e_{#4}e_{#5}]\ldots] }

\newcommand{\eiajthree}[3]{[e_{#1}[\ldots[e_1[e_1\ldots[e_{#2}e_{#3}]\ldots] }
\newcommand{\one}{[e_1[e_1[e_2e_3]]]}

\newcommand{\depi}[2]{[e_{\bar{1}}[e_2[\ldots[e_{#1}e_{#2}]\ldots]}
\newcommand{\dei}[2]{[e_{1}[e_2[\ldots[e_{#1}e_{#2}]\ldots]}
\newcommand{\dejj}[2]{[e_{\bar{1}}[e_1[\ldots[e_{#1}e_{#2}]\ldots]}
\newcommand{\deij}[3]{[e_{#1}[\ldots[e_{#2}e_{#3}]\ldots]}
\newcommand{\deiajthree}[3]{[e_{#1}[\ldots[e_{\bar{1}}[e_1\ldots[e_{#2}e_{#3}]\ldots] }
\newcommand{\deiajthreec}[3]{[f_{#1}[\ldots[f_{1}[f_1\ldots[f_{#2}f_{#3}]\ldots] }
\newcommand{\deiajthreet}[3]{[e_{#1}[\ldots[\frac{e_{1}+e_{\bar{1}}}{2}[\frac{e_1+e_{\bar{1}}}{2}\ldots[e_{#2}e_{#3}]\ldots] }
\newcommand{\deiajthreetm}[3]{[e_{#1}[\ldots[e_{1}-e_{\bar{1}}[e_1-e_{\bar{1}}[\ldots[e_{#2}e_{#3}]\ldots] }
\newcommand{\dop}{e_{\bar{1}}}
\newcommand{\don}{e_{1}}

\newcommand{\feij}[3]{[e_{#1}[\ldots[e_{#2}e_{#3}]\ldots]}

\newcommand{\feiajthree}[3]{[e_{#1}[\ldots[e_1\ldots[e_{#2}e_{#3}]\ldots] }
\newcommand{\fett}{[e_2[e_1e_2]]}
\newcommand{\fetth}{[e_2[e_1[e_2e_3]]]}
\newcommand{\fetf}{[e_2[e_1[e_2[e_3e_4]]]]}
\newcommand{\fethth}{[e_3[e_2[e_1[e_2e_3]]]]}
\newcommand{\fethf}{[e_3[e_2[e_1[e_2[e_3e_4]\ldots]}
\newcommand{\feff}{[e_4[e_3[e_2[e_1[e_2[e_3e_4]\ldots]}

\renewcommand{\Im}{\mathop{\mathrm{Im}}}

\newcommand{\cB}{\mathcal{B}}

\newcommand{\cD}{\mathcal{D}}

\newcommand{\fe}{\mathfrak{e}}

\author{YI SU}
\title{Electrical Lie Algebra of Classical Types}

\begin{document}

\maketitle

\begin{abstract}
We investigate the structure of electrical Lie algebras of finite Dynkin type. These Lie algebras were introduced by Lam-Pylyavskyy in the study of \textit{circular planar electrical networks}. The corresponding Lie group acts on such networks via some combinatorial operations studied by Curtis-Ingerman-Morrow and Colin de Verdi\`{e}re-Gitler-Vertigan. Lam-Pylyavskyy studied the electrical Lie algebra of type $A$ of even rank in detail, and gave a conjecture for the dimension of electrical Lie algebras of finite Dynkin types. We prove this conjecture for all classical Dynkin types, that is, $A$, $B$, $C$, and $D$.  Furthermore, we are able to explicitly describe the structure of the corresponding electrical Lie algebras as the semisimple product of the symplectic Lie algebra with its finite dimensional irreducible representations.
\end{abstract}

\section{Introduction}

The study of electrical networks dated back to Georg Ohm and Gustav Kirchhoff more than a century ago, and it is still a classical object in the study of many branches of mathematics including graph theory (see for example \cite{kw}). It also has many applications in other fields including material science and medical imaging (see for example \cite{bvm}). \vspace{0.1in}

In this paper, our main object of study, \emph{the electrical Lie algebra}, originates from the study of \textit{circular planar electrical networks}.

\begin{center}
\begin{tikzpicture}[scale=0.6]

\node at (0:4.6) {$\bar{5}$};
\node at (45:4.6) {$\bar{4}$};
\node at (90:4.6) {$\bar{3}$};
\node at (135:4.6) {$\bar{2}$};
\node at (180:4.6) {$\bar{1}$};
\node at (225:4.6) {$\bar{8}$};
\node at (270:4.6) {$\bar{7}$};
\node at (315:4.6) {$\bar{6}$};

\draw[dashed] (0,0) circle (4cm);
\draw (315:4) -- node [above] {$3$} (270:2);
\draw (270:2) -- node [left] {$1$} (270:4);
\draw (270:2) -- node [above] {$2$} (0:4);
\draw (150:2) -- node [left] {$2$} (225:4);
\draw (150:2) -- node [above] {$5$} (90:4);
\draw (150:2) -- node [above] {$1$} (270:2);
\draw (0:4) -- node [left] {$4$} (45:4);
\draw (135:4) -- node [right] {$1$} (180:4);
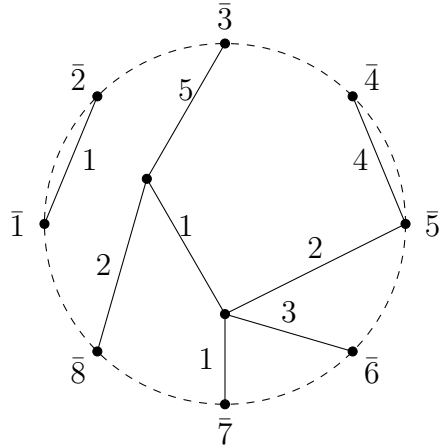
\draw (270:5.6)node{\text{Circular Planar Electrical Network}};
\filldraw[black] (0:4) circle (0.1cm);
\filldraw[black] (45:4) circle (0.1cm);
\filldraw[black] (90:4) circle (0.1cm);
\filldraw[black] (135:4) circle (0.1cm);
\filldraw[black] (180:4) circle (0.1cm);
\filldraw[black] (225:4) circle (0.1cm);
\filldraw[black] (270:4) circle (0.1cm);
\filldraw[black] (315:4) circle (0.1cm);
\filldraw[black] (150:2) circle (0.1cm);
\filldraw[black] (270:2) circle (0.1cm);

;

\end{tikzpicture}
\end{center}
The numbers on the edges are conductances. The vertices on the outer circular boundary are ordered, and called boundary vertices, and the rest are called interior vertices. Edges cannot cross with other edges. When voltage is put on the boundary vertices, electrical current will flow in or out of the boundary vertices.\vspace{0.1in}

Curtis-Ingerman-Morrow \cite{cim} and Colin de Verdi\`{e}re-Gitler-Vertigan \cite{dvgv} gave a robust theory of circular planar electrical networks. In \cite{cim} and \cite{dvgv}, two operations of \textit{adjoining a boundary spike} and \textit{adjoining a boundary edge} to a circular planar electrical network were studied:

\begin{center}
\begin{tikzpicture}
\draw[dashed] (0,0) circle (1.5cm);
\draw (45:1.5) -- (45:2.2);
\filldraw[black] (0:1.5) circle (0.07cm);
\filldraw[black] (45:1.5) circle (0.07cm);
\filldraw[black] (90:1.5) circle (0.07cm);
\filldraw[black] (135:1.5) circle (0.07cm);
\filldraw[black] (180:1.5) circle (0.07cm);
\filldraw[black] (225:1.5) circle (0.07cm);
\filldraw[black] (270:1.5) circle (0.07cm);
\filldraw[black] (315:1.5) circle (0.07cm);
\filldraw[black] (45:2.2) circle (0.07cm);
\draw (270:2.2) node{Adjoining a boundary spike};

\begin{scope}[shift={(6,0)}]
\draw[dashed] (0,0) circle (1.5cm);
\draw (45:1.5) -- (0:1.5);
\filldraw[black] (0:1.5) circle (0.07cm);
\filldraw[black] (45:1.5) circle (0.07cm);
\filldraw[black] (90:1.5) circle (0.07cm);
\filldraw[black] (135:1.5) circle (0.07cm);
\filldraw[black] (180:1.5) circle (0.07cm);
\filldraw[black] (225:1.5) circle (0.07cm);
\filldraw[black] (270:1.5) circle (0.07cm);
\filldraw[black] (315:1.5) circle (0.07cm);
\draw (270:2.2) node{Adjoining a boundary edge};
\end{scope}

\end{tikzpicture}
\end{center}

These two operations generate the set of circular planar electrical networks modulo the electrical equivalences. Lam-Pylyavskyy \cite{lp} viewed these actions as one parameter subgroups of a Lie group action, namely the electrical Lie group (of type $A$). Then they define electrical Lie algebras of finite Dynkin type:

\begin{definition}
Let $X$ be a Dynkin diagram of finite type, $I=I(X)$ be the set of nodes in $X$, and $A=(a_{ij})$ be the associated Cartan matrix . Define the electrical Lie algebra $\fe_{X}$ associated to $X$ to be the Lie algebra generated by $\{e_i\}_{i\in I}$ modulo the relations 

\[\mathfrak{ad}(e_i)^{1-a_{ij}}(e_j)=
\begin{cases}0  &\text{ if } a_{ij}\neq -1,\\
                    2e_i&\text{ if }a_{ij}=-1.
\end{cases}\]
where $\mathfrak{ad}$ is the adjoint representation.
\end{definition}

Note that our convention for $a_{ij}$ is that if the simple root corresponding to $i$ is shorter than the one corresponding to an adjacent node $j$, then $|a_{ij}|>1$. Equivalently, the arrows in the  Dynkin diagram point towards the nodes which correspond to the shorter roots.\vspace{0.1in}

These relations can be seen as a deformation of the upper half of semisimple Lie algebras. For ordinary semisimple Lie algebras,  the corresponding relations are

\[\mathfrak{ad}(e_i)^{1-a_{ij}}(e_j)=0\ \ \ \forall i,j.\]\vspace{0.05in}

\begin{center}
\begin{tikzpicture}[scale=0.7]
\filldraw[black] (0:0) circle (0.1cm);
\filldraw[black] (90:2.5) circle (0.1cm);
\filldraw[black] (210:2.5) circle (0.1cm);
\filldraw[black] (330:2.5) circle (0.1cm);
\draw(0:0)--(90:2.5);
\draw(0:0)--(210:2.5);
\draw(0:0)--(330:2.5);

\begin{scope}[shift={(3.5,0)}]
\node at (0:0) {$\longleftrightarrow$};
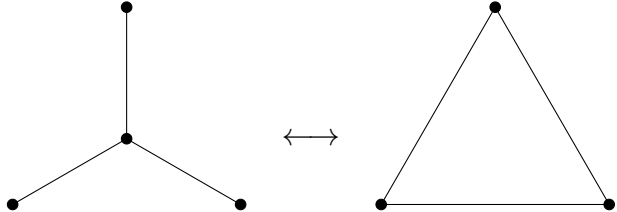
\node at (270:2.5) {Star-Triangle Transformation};
\end{scope}

\begin{scope}[shift={(7,0)}]
\filldraw[black] (90:2.5) circle (0.1cm);
\filldraw[black] (210:2.5) circle (0.1cm);
\filldraw[black] (330:2.5) circle (0.1cm);
\draw(90:2.5)--(210:2.5);
\draw(90:2.5)--(330:2.5);
\draw(210:2.5)--(330:2.5);
\end{scope}
\end{tikzpicture}
\end{center}

In the case of the electrical Lie algebra of type $A$, the famous star-triangle (or Yang-Baxter) transformation of electrical networks translates into ``\textit{the electrical Serre relation}'': 
\[[e_i,[e_i,e_{i\pm1}]]=-2e_i,\]
whereas the usual Serre relation for the semisimple Lie algebra of type $A$ is
\[[e_i,[e_i,e_{i\pm1}]]=0.\]
where $i\in [n]$ are labels of the nodes of the Dynkin diagram $A_n$. Lam-Pylyavskyy looked at the algebraic structure of electrical Lie groups and Lie algebras of finite Dynkin type. They showed that $\fe_{A_{2n}}$ is semisimple and isomorphic to the symplectic Lie algebra $\mathfrak{sp}_{2n}$. Moreover, they conjectured that the dimension of the electrical Lie algebra $\fe_{X}$ equals the number of positive roots $|\Phi(X)^+|$, where $\Phi(X)^+$ is the set of positive roots of root system $\Phi(X)$ with Dynkin diagram $X$. \vspace{0.1in}

In this current paper, we will not only prove Lam-Pylyavskyy's conjecture regarding the dimension for all classical types, but also explore the structure of certain electrical Lie algebras of classical types.\vspace{0.1in}

In fact, the semisimplicity of $\fe_{A_{2n}}$ is not a general property of electrical Lie algebras. For example, in this paper we will see that $\fe_{C_{2n}}$ has a nontrivial solvable ideal. This also makes the structure of such Lie algebras difficult to describe. There is no uniform theory for electrical Lie algebras of classical type, so we will explore the structure of such electrical Lie algebras in a case by case basis:\vspace{0.1in}

Consider three irreducible representations of $\mathfrak{sp}_{2n}$: let $V_{\mathbf{0}}$ be the trivial representation, $V_{\nu}$ be the standard representation, that is, the highest weight vector $\nu=\omega_1$, and $V_{\lambda}$ be the irreducible representation with the highest weight vector $\lambda=\omega_1+\omega_2$, where $\omega_1$, and $\omega_2$ are the first and second fundamental weights.\vspace{0.1in}

For type $A$, the electrical Lie algebra of even rank $\fe_{A_{2n}}$ is isomorphic to $\mathfrak{sp}_{2n}$ \cite{lp}. We show that $\fe_{A_{2n+1}}$ is isomorphic to an extension of $\mathfrak{sp}_{2n}\ltimes V_{\nu}$ by the $\mathfrak{sp}_{2n}$-representation $V_{\mathbf{0}}$. Note that $\fe_{A_{2n+1}}$ is also isomorphic to the odd symplectic Lie algebra $\mathfrak{sp}_{2n+1}$ studied by Gelfand-Zelevinsky \cite{gz} and Proctor \cite{rp}.\vspace{0.1in}

For type $B$, we show that $\fe_{B_{n}}\cong\mathfrak{sp}_{n}\oplus \mathfrak{sp}_{n-1}$ by constructing an isomorphism between these two Lie algebras, where the odd symplectic Lie algebra is the same as the one appearing in the case of type $A$.\vspace{0.1in}

For type $C$, we first consider the case of even rank. We find an abelian ideal $I\subset \fe_{C_{2n}}$ and prove that this quotient $\fe_{C_{2n}}/I$ is isomorphic to $\fe_{A_{2n}}$. So we can define a Lie group action of $\fe_{A_{2n}}$ (or $\mathfrak{sp}_{2n}$) on $I$. Consequently, we show that $\fe_{C_{2n}}$ is isomorphic to $\mathfrak{sp}_{2n}\ltimes (V_{\lambda}\oplus V_{\mathbf0})$. As for the odd case, it is a Lie subalgebra of $\fe_{C_{2n+2}}$, so we are able to conclude that its dimension is $(2n+1)^2$, the number of positive roots $|\Phi^+(C_{2n+1})|$.\vspace{0.1in}

For type $D$, we find that $\fe_{D_{n+1}}$ contains a Lie subalgebra isomorphic to $\fe_{C_{n}}$, and use the structure theorem of $\fe_{C_{2n}}$ to find that the dimension of $\fe_{D_{2n+1}}$ is equal to the number of positive roots $|\Phi^+(D_{2n+1})|$. Similarly to type $C$, we are also able to conclude that the dimension of $\fe_{D_{2n}}$ is the one expected as in the conjecture.\vspace{0.1in}

In a similar manner to the connection between \textit{circular planar electrical networks} and the electrical Lie algebra $\fe_{A_{n}}$, we find that $\fe_{B_{n}}$ also has strong connection with \textit{mirror symmetric circular planar eletrical networks} via similar combinatorial operations on such networks. We will only focus on the type $B$ electrical Lie algebra in the current paper, and introduce this analogue and investigate the properties of mirror symmetric circular planar electrical networks in a future paper.\vspace{0.1in}

The structure of this paper goes as follows: Section 2 is  the outline of the proofs of structure theorems of type $A$, $B$, $C$, and $D$ electrical Lie algebras, whereas the proofs of some technical lemmas in Subsection \ref{typec} and \ref{typed} are left in Section 3.\vspace{0.3in}

\section{Electrical Lie Algebra of Classical Types}
\label{classicaltypes}

\subsection{Preliminary Lemmas and Propositions}
Lam-Pylyavskyy \cite{lp} prove the following proposition that gives an upper bound for the dimension of electrical Lie algebras of finite Dynkin type.

\begin{proposition}[\cite{lp}]
\label{dimupperbound}
 Let $X$ be a Dynkin diagram of finite type, $\Phi^+(X)$ be the set of positive roots of $X$. The dimension of $\fe_X$ is less than or equal to $|\Phi^+(X)|$. Moreover, $\fe_{X}$ has a spanning set indexed by positive roots in $|\Phi^+(X)|$. 
\end{proposition}

\begin{remark} 
\label{spanningset}
In the proof of Proposition \ref{dimupperbound}, Lam-Pylyavskyy gave a concrete spanning set of this Lie algebra: Let $\aph$ be any element in $\Phi^+(X)$. Write $\aph=\aph_{i_1}+\aph_{i_2}+\ldots+\aph_{i_t}$ such that each $\aph_{i_j}$ is a simple root, and $\sum_{j=1}^{s}\aph_{i_j}$ is a positive root for all $s\in [t]$. Set $e_{\aph}:=[e_{i_1}[e_{i_2}[\ldots[e_{i_{t-1}}e_{i_t}]\ldots]$. Then the set $\{e_{\aph}\}_{\aph\in \Phi^+(X)}$ is a spanning set of $\fe_{X}$.
\end{remark}

We will also need the following lemma in exploring the structure of electrical Lie algebras:

\begin{lemma}\label{semisimpleaction}
Let $L$ be a Lie algebra, and $I$  an ideal of $L$ with $[I,I]=0$. Then the quotient Lie algebra $L/I$ has a Lie algebra action on $I$.
\end{lemma}

\proof

Let $\bar{a}\in L/I$, where $a\in L$. Let $x\in I$. Define $\bar{a}\cdot x=[a,x]$. Let $b\in L$ such that $\bar{b}=\bar{a}\in L/I$. Hence, $a-b=y\in I$. $(\bar{a}-\bar{b})\cdot x=[a-b,x]=[y,x]=0$ because of $[I,I]=0$, which shows the well-definedness. Since this action is induced by the adjoint representation of $L$ on $I$, it defines a valid representation of $L/I$. $\square$\vspace{0.1in}

We will show that the upper bound in Proposition \ref{dimupperbound} is the correct dimension for each Dynkin diagram of classical type. Furthermore, we will also give the explicit structure of electrical Lie algebras $\fe_{A_n}$, $\fe_{B_{n}}$, and $\fe_{C_{2n}}$. The notations will only be used within each following subsection of Section \ref{classicaltypes}.\vspace{0.1in}

\subsection{Type $A$}
\label{typea}
\vspace{0.1in}

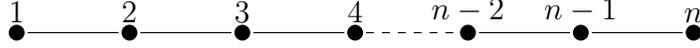
\begin{figure}
\begin{tikzpicture}[scale=1.5]
\fill

(1,0) node(1){} circle(2pt)
(2,0) node(2){}circle(2pt)
(3,0) node(3){} circle(2pt)
(4,0) node(4){} circle(2pt)
(5,0) node(5){} circle(2pt)
(6,0) node(6){} circle(2pt)
(7,0) node(7){} circle(2pt)
(1,0) node(n1)[above]{1}
(2,0) node(n2)[above]{2}
(3,0) node(n3)[above]{3}
(4,0) node(n4)[above]{4}
(5,0) node(n5)[above]{$n-2$}
(6,0) node(n6)[above]{$n-1$}
(7,0) node(n7)[above]{$n$}
;

\draw (1) -- (2);
\draw (2) -- (3);
\draw (3) -- (4);
\draw[dashed] (4) -- (5);
\draw (5) -- (6);
\draw (6) -- (7);

\end{tikzpicture}
\caption{Dynkin Diagram of $A_{n}$}
\label{typeadynkin}
\end{figure}

By definition, Lie algebra $\fe_{A_{n}}$ is generated by generators $\{e_i\}_{i=1}^n$ under the relations (Figure \ref{typeadynkin})
\begin{alignat*}{3}
\label{relationa}
[e_i,e_j]        &=0          &&\text{ if } |i-j|\ge 2,\\
[e_i,[e_i,e_j]]&=-2e_i   &&\text{ if } |i-j|=1.
\end{alignat*}

\vspace{0.1in}

Lam-Pylyavskyy studied the structure of $\fe_{A_{2n}}$:

\begin{theorem}[\cite{lp}]
\label{ea2n}
We have 
\[\fe_{A_{2n}}\cong \mathfrak{sp}_{2n}.\]
\end{theorem}\vspace{0.1in}

We will explore the structure of $\fe_{A_{2n+1}}$:\vspace{0.1in}

\begin{proposition}
\label{oddtypeadim}
The dimension of $\fe_{A_{2n+1}}$ is $(n+1)(2n+1)$.
\end{proposition} 

\proof

Let $\{e_i\}_{i=1}^{2n+1}$ and $\{\tilde{e}_{i}\}_{i=1}^{2n+2}$ be generators of $\fe_{A_{2n+1}}$ and $\fe_{A_{2n+2}}$ respectively. Then there is a Lie algebra homomorphism:
\[\psi:\fe_{A_{2n+1}}\longrightarrow\fe_{A_{2n+2}},\quad e_i\longmapsto \tilde{e}_i.\] 

We claim this is an injection.\vspace{0.1in}

By Remark \ref{spanningset}, $\{e_{\aph}\}_{\aph\in \Phi^+(A_{2n+1})}$ and $\{\tilde{e}_{\aph}\}_{\aph\in \Phi^+(A_{2n+2})}$ are spanning sets of $\fe_{A_{2n+1}}$ and $\fe_{A_{2n+2}}$ respectively. By Theorem \ref{ea2n}, we know that $\{\tilde{e}_{\aph}\}_{\aph\in\Phi^+(A_{2n+2})}$ is a basis of $\fe_{A_{2n+2}}$. On the other hand $\psi(e_{\aph})=\tilde{e}_{\aph}$, so $\{e_{\aph}\}_{\aph\in\Phi^+(A_{2n+1})}$ is a basis of $\fe_{A_{2n+1}}$. Therefore, $\psi$ is an injective Lie algebra homomorphism, and the dimension of $\fe_{A_{2n+1}}$ is $(n+1)(2n+1)$. $\square$\vspace{0.2in}

We consider a Lie subalgebra of $\mathfrak{sp}_{2n+2}$ introduced by Gelfand-Zelevinsky \cite{gz}. The definition of this Lie algebra that we will use comes from Proctor \cite{rp}. Let $V=\mathbb{C}^{2n+2}$ be a complex vector space with standard basis $\{x_i\}_{i=1}^{2n+2}$. Let $\{y_i\}_{i=1}^{2n+2}$ be the corresponding dual basis of $V^*$. Also let $B$ be a nondegenerate skew-symmetric bilinear form of $V$, and $\beta$ be one nonzero element in $V^*$. Let $GZ(V,\beta,B)$ be the Lie subgroup of $GL(V)$ which preserves both $\beta$ and $B$. Its Lie algebra is denoted as $gz(V,\beta,B)$. Since $Sp_{2n+2}=\{M\in GL(V)| M \text{ preserves } B\}$, we have $gz(V,\beta,B)\subset \mathfrak{sp}_{2n+2}$. Now we fix $\beta=y_{n+1}$ and the matrix representing $B$ to be $\begin{pmatrix}
0 & I_{n+1}\\
-I_{n+1} & 0
\end{pmatrix}$. Then we define the \emph{odd symplectic Lie algebra} $\mathfrak{sp}_{2n+1}$ to be $gz(\mathbb{C}^{2n+2},y_{n+1}, B)$.\vspace{0.2in}

\begin{theorem}
\label{oddtypeastructure1}

We have $\fe_{A_{2n+1}}\cong \mathfrak{sp}_{2n+1}$ as Lie algebras.
\end{theorem}

\proof We use a specific  isomorphism of $\fe_{A_{2n+2}}$ and $\mathfrak{sp}_{2n+2}$ in Theorem 3.1 of \cite{lp}. Let $\epsilon_{i}\in \mathbb{C}^{2n+2}$ denote the  column vector with 1 in the $i$th position and 0 elsewhere. Let $a_1=\epsilon_1,a_2=\epsilon_1+\epsilon_2,\ldots,a_{n+1}=\epsilon_{n}+\epsilon_{n+1}$, and $b_1=\epsilon_1,b_2=\epsilon_2,\ldots,b_{n+1}=\epsilon_{n+1}$.  Now define $\phi:\fe_{A_{2n+2}}\longrightarrow \mathfrak{sp}_{2n+2}$ as follows:

\[e_{2i-1}\mapsto\begin{pmatrix}0& a_i\cdot a_i^T\\ 0& 0\end{pmatrix}, \qquad e_{2i}\mapsto\begin{pmatrix} 0 & 0 \\ b_i\cdot b_i^T & 0\end{pmatrix}.\]
and extend this to a Lie algebra homomorphism.\vspace{0.1in}

It can be found that $\mathfrak{sp}_{2n+1}$ consists of the matrices of $\mathfrak{sp}_{2n+2}$ whose entries in $(n+1)$-th column and $(2n)$-th row are all zero. The dimension of this Lie algebra is $(n+1)(2n+1)$.\vspace{0.1in}

Due to Theorem \ref{oddtypeadim}, both $\fe_{A_{2n+1}}$ and $\mathfrak{sp}_{2n+1}$ have dimension $(n+1)(2n+1)$. Since $\phi\circ\psi$ is injective, it suffices to show $\phi\circ\psi(\fe_{A_{2n+1}})\subset \mathfrak{sp}_{2n+1}$. Because $\phi\circ\psi(\fe_{A_{2n+1}})\subset \mathfrak{sp}_{2n+2}$ and \[\mathfrak{sp}_{2n+2}=\left\{\begin{pmatrix}
m & n\\
p & q
\end{pmatrix}\Big| m,n,p,q \text{ are } (n+1)\times(n+1)\text{ matrices}, m=-q^T, p=p^T, n=n^T\right\},\]
we only need to show the entries of the $(n+1)$-th column of $\phi\circ\psi(\fe_{A_{2n+1}})$ are all zeros. Clearly the entries of the $(n+1)$-th column of $\phi\circ\psi(e_{i})$ for $i=1,2,\ldots, 2n+1$ are all zeros. Now assume $
R_i=\begin{pmatrix}A_i & B_i \\C_i & D_i\end{pmatrix}$ for $i=1,2$ be block matrices, where all entries are the matrices of size $(n+1)\times(n+1)$, and the entries of $(n+1)$-th column of $A_i$ and $C_i$ are all zeros. Notice that if $M,N$ are two square matrices, and the entries of last column of $N$ are all zeros, then the last column of $MN$ is a zero column vector. Now
\[R_1R_2=\begin{pmatrix}A_1 & B_1 \\ C_1 & D_1\end{pmatrix}\begin{pmatrix}A_2 & B_2 \\ C_2 & D_2\end{pmatrix}
=\begin{pmatrix}A_1A_2+B_1C_2 & A_1B_2+B_1D_2\\ A_2C_1+D_1C_2 & B_2C_1+D_1D_2\end{pmatrix}.\] 
Therefore, the last columns of $A_1A_2+B_1C_2$ and $A_2C_1+D_1C_2$ are both zero. So the $(n+1)$-th column of the product $R_1R_2$ is also zero. Since this property of having $(n+1)$-th column being zero is closed among set of $(2n+2)\times(2n+2)$ matrices under matrix multiplication, it is also closed under Lie bracket. Thus, $\phi\circ\psi(\fe_{A_{2n+1}})\subset \mathfrak{sp}_{2n+1}$. $\square$\vspace{0.1in}

\begin{theorem}
\label{oddtypeastructure2}

We have that $\fe_{A_{2n+1}}$ is isomorphic to an extension of $\mathfrak{sp}_{2n}\ltimes V_{\nu}$ by the representation $V_{\mathbf{0}}$, where $V_{\nu}$ and $V_{\mathbf{0}}$ are standard and trivial representation of $\mathfrak{sp}_{2n}$ respectively. In another word, there is a short exact sequence 
\[0\longrightarrow \mathbb{C}\longrightarrow \fe_{A_{2n+1}}\longrightarrow \mathfrak{sp}_{2n}\ltimes \mathbb{C}^{2n}\longrightarrow 0.\]

\end{theorem}

\proof We will use the matrix presentation of $\fe_{A_{2n+1}}$ in Theorem \ref{oddtypeastructure1} to prove this theorem.

Let \[A=\begin{pmatrix}
0          & \cdots    & 0                & 0\\ 
\vdots  & \ddots   & \vdots        & \vdots\\
0          & \cdots    & 0                & 0\\
a_1 & \cdots  & a_{n}  & 0
\end{pmatrix},\quad
B=\begin{pmatrix}
0          & \cdots    & 0                      & b_{1}\\ 
\vdots  & \ddots   & \vdots              & \vdots\\
0          & \cdots    & 0                      & b_{n}\\
b_1 & \cdots  & b_n  & b_{n+1}
\end{pmatrix},\quad C=0.\] \vspace{0.1in}

Let $I$ be the set of matrices of the form 
$\begin{pmatrix}A & B\\
C & -A^t
\end{pmatrix}$, where $a_i,b_i$'s are arbitrary. It is clear that $I$ is an ideal of $\fe_{A_{2n+1}}$, and $\fe_{A_{2n+1}}/I\cong \fe_{A_{2n}}$. On the other hand, let $I'$ be the one dimensional subspace of $\fe_{A_{2n+1}}$ generated by 
$\begin{pmatrix}
0 & E_{n+1,n+1}\\
0 & 0
\end{pmatrix}$. Then $I'$ is in the center of $\fe_{A_{2n+1}}$. With calculation we know that $[I,I]\subset I'$. Let $\bar{I}$ be the image of $I$ in $\fe_{A_{2n+1}}/I'$. Thus in the quotient algebra $\fe_{A_{2n+1}}/I'$, we have $[\bar{I},\bar{I}]=0$. By Lemma \ref{semisimpleaction}, this shows that $\fe_{A_{2n+1}}/I'/\bar{I}\cong \fe_{A_{2n}}\cong \mathfrak{sp}_{2n}$ has an action on $\bar{I}$. Since $\bar{I}$ is cyclic and $\dim \bar{I}=2n$, it has to be isomorphic to the standard representation $V_{\nu}$. Thus, 
\[\fe_{A_{2n+1}}/I'\cong \spn_{2n}\ltimes V_{\nu}.\]

Because $I'$ is in the center of $\fe_{A_{2n+1}}$, it is isomorphic to the trivial representation $V_{\mathbf{0}}$ of $\spn_{2n}$.\vspace{0.1in} 

Therefore, $\fe_{A_{2n+1}}$ is isomorphic to an extension of $\mathfrak{sp}_{2n}\ltimes V_{\nu}$ by the trivial representation $V_{\mathbf{0}}$. According to the above matrix presentation of $\fe_{A_{2n+1}}$, we can find a short exact sequence 

\[0\longrightarrow \mathbb{C}\longrightarrow \fe_{A_{2n+1}}\longrightarrow \mathfrak{sp}_{2n}\ltimes \mathbb{C}^{2n}\longrightarrow 0.\] $\square$\vspace{0.3in}

This concludes the study of the electrical Lie algebra of type $A$.\vspace{0.3in}

\subsection{Type $B$}
\label{typeb}

The electrical Lie algebra $\fe_{B_{n}}$ is generated by $\{e_1,e_2,\ldots,e_n\}$, with relations (Figure \ref{typebdynkin}):

\begin{alignat*}{3}
&[e_i,e_j]=0 &&\text{ if } |i-j|\ge 2\\
&[e_i,[e_i,e_j]]=-2e_i &&\text{ if } |i-j|=1, i\neq2 \text{ and } j\neq 1\\
&[e_2,[e_2,[e_2,e_1]]]=0\\
\end{alignat*}

\begin{figure}
\begin{tikzpicture}[scale=1.5]
\fill

(1,0) node(1){} circle(2pt)
(2,0) node(2){}circle(2pt)
(3,0) node(3){} circle(2pt)
(4,0) node(4){} circle(2pt)
(5,0) node(5){} circle(2pt)
(6,0) node(6){} circle(2pt)
(7,0) node(7){} circle(2pt)
(1,0) node(n1)[above]{1}
(2,0) node(n2)[above]{2}
(3,0) node(n3)[above]{3}
(4,0) node(n4)[above]{4}
(5,0) node(n5)[above]{$n-2$}
(6,0) node(n6)[above]{$n-1$}
(7,0) node(n7)[above]{$n$}
;

\path[doublearrow] (2) to (1);
\draw (2) -- (3);
\draw (3) -- (4);
\draw[dashed] (4) -- (5);
\draw (5) -- (6);
\draw (6) -- (7);

\end{tikzpicture}
\caption{Dynkin Diagram of $\cB_{n}$}
\label{typebdynkin}
\end{figure}

Let $\{f_1,f_2,\ldots,f_n\}$ be a generating set of $\fe_{A_{n}}$. Then $\{f_2,\ldots,f_n\}$ is a generating set of $\fe_{A_{n-1}}$. Now consider a new Lie algebra $\fe_{A_{n}}\oplus \fe_{A_{n-1}}$: the underlying set is $(a,b)\in \fe_{A_{n}}\times \fe_{A_{n-1}}$, and the Lie bracket operation is $[(a,b),(c,d)]=([a,c],[b,d])$.\vspace{0.2in}

Define a map $\phi: \fe_{B_{n}}\longrightarrow \fe_{A_{n}}\oplus\fe_{A_{n-1}}$ as follows:
\begin{align*}
&\phi(e_1)=(f_1,0),\qquad\phi(e_k)=(f_k,f_k)\ \ \forall k\ge 2\\
&\phi([e_{i_1}[e_{i_2}[\ldots[e_{i_{s-1}}e_{i_{s}}]\ldots])=[\phi(e_{i_1})[\phi(e_{i_2})[\ldots[\phi(e_{i_{s-1}})\phi(e_{i_{s}})]\ldots]
\end{align*}\vspace{0.01in}

\begin{theorem}
\label{typebstructure}

$\phi$ is a Lie algebra isomorphism. Therefore, we have
\[\fe_{B_{n}}\cong\spn_{n}\oplus\spn_{n-1},\]
where the odd symplectic Lie algebra is defined in Section \ref{typea}.
\end{theorem}

\proof

First of all we would like to prove $\phi$ is a Lie algebra homomorphism. It suffices to show that $\phi(e_k)$'s also satisfy the defining relation of $\fe_{B_{n}}$.

For $k\ge 2$, we have 
\begin{align*}
[\phi(e_k)[\phi(e_k)\phi(e_{k+1})]]=&[(f_k,f_k),[(f_k,f_k),(f_{k+1},f_{k+1})]]\\
                                                         =&([f_k[f_kf_{k+1}]],[f_k[f_kf_{k+1}]]\\
                                                         =&-2(f_k,f_k)=-2\phi(e_k).
\end{align*}

Similarly, for $k\ge 3$
\[[\phi(e_k)[\phi(e_k)\phi(e_{k-1})]=-2\phi(e_{k-1}).\]

And
\begin{align*}
[\phi(e_1)[\phi(e_1)\phi(e_2)]]=&[(f_1,0),[(f_1,0),(f_2,f_2)]]\\
                                                 =&([f_1[f_1f_2]],0)=-2(f_1,0)\\
                                                 =&-2\phi(e_1),
\end{align*}
\begin{align*}
[\phi(e_2)[\phi(e_2)[\phi(e_2)\phi(e_1)]]]=&[(f_2,f_2),[(f_2,f_2),[(f_2,f_2),(f_1,0)]]]\\
                                                                  =&[(f_2,f_2),(-2f_2,0)]=0.
\end{align*}

It is also clear that if $|i-j|\ge2$, $[\phi(e_i)\phi(e_j)]=0$. Therefore, this is a Lie algebra homomorphism.\vspace{0.1in}

Next we claim that $\phi$ is surjective. We already know that $(f_1,0)\in \Im(\phi)$, so it suffices to show that $(f_k,0)$, $(0,f_k)\in \Im(\phi)$ for all $k\ge 2$. We go by induction. Base case:
\begin{align*}
&\phi(-\dfrac{1}{2}[e_2[e_2e_1]])=-\dfrac{1}{2}[(f_2,f_2),[(f_2,f_2),(f_1,0)]]=(f_2,0),\\
&\phi(e_2+\dfrac{1}{2}[e_2[e_2e_1]])=(f_2,f_2)-(f_2,0)=(0,f_2).
\end{align*}
So $(f_2,0),\ (0,f_2)\in\Im(\phi)$.\vspace{0.1in}

Now assume this is true for $k$. Without loss of generality, say $\phi(y)=(f_k,0)$. Then
\begin{align*}
&\phi(-\dfrac{1}{2}[e_{k+1}[e_{k+1}y]])=-\dfrac{1}{2}[(f_{k+1},f_{k+1}),[(f_{k+1},f_{k+1}),(f_k,0)]]=(f_{k+1},0),\\
&\phi(e_{k+1}+\dfrac{1}{2}[e_{k+1}[e_{k+1}y]])=(f_{k+1},f_{k+1})-(f_{k+1},0)=(0,f_{k+1}).
\end{align*}

Thus, the claim is true. Note that $\dim \fe_{A_{n}}\oplus\fe_{A_{n-1}}=\binom{n+1}{2}+\binom{n}{2}=n^2$. By Proposition \ref{dimupperbound}, $\dim \fe_{B_{n}}\le n^2$. Then we get
\[n^2=\dim \fe_{A_{n}}\oplus \fe_{A_{n-1}}\le \dim \fe_{B_{n}}\le n^2.\]
So we achieve equality, and $\phi$ is an isomorphism. By Section \ref{typea}, we know that $\fe_{A_{n}}\cong \spn_{n}$ for all $n$. Hence,
\[\fe_{B_{n}}\cong\spn_{n}\oplus\spn_{n-1}.\]  $\square$\vspace{0.3in}

\subsection{Type $C$}
\label{typec}

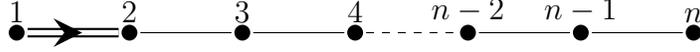
\begin{figure}
\begin{tikzpicture}[scale=1.5]
\fill

(1,0) node(1){} circle(2pt)
(2,0) node(2){}circle(2pt)
(3,0) node(3){} circle(2pt)
(4,0) node(4){} circle(2pt)
(5,0) node(5){} circle(2pt)
(6,0) node(6){} circle(2pt)
(7,0) node(7){} circle(2pt)
(1,0) node(n1)[above]{1}
(2,0) node(n2)[above]{2}
(3,0) node(n3)[above]{3}
(4,0) node(n4)[above]{4}
(5,0) node(n5)[above]{$n-2$}
(6,0) node(n6)[above]{$n-1$}
(7,0) node(n7)[above]{$n$}
;

\path[doublearrow] (1) to (2);
\draw (2) -- (3);
\draw (3) -- (4);
\draw[dashed] (4) -- (5);
\draw (5) -- (6);
\draw (6) -- (7);

\end{tikzpicture}
\caption{Dynkin Diagram of $C_{n}$}
\label{typecdynkin}
\end{figure}

The electrical Lie algebra of type $C_n$ is generated by generators $\{e_i\}_{i=1}^n$ with relations (Figure \ref{typecdynkin}):
\begin{alignat*}{4}
&[e_i,e_j]                     =0,        &&\text{ if } |i-j|\ge 2,\\
&[e_i,[e_i,e_j]]            =-2e_i,   &&\text{ if } |i-j|=1,\ i\neq 1,\\
&[e_1[e_1[e_1e_2]]]   =0.     && \text{ }
\end{alignat*}

By Remark \ref{spanningset}, there is a spanning set of $\fe_{C_{n}}$ indexed by the positive roots $\Phi^+(C_{n})$. More precisely, the spanning set is 
\[\{\eiaj{i}{j}:1\leq i<j\leq n\}\cup \{\eij{i}{j-1}{j}:1\leq i\leq j\leq n\}.\] 
The way in which the generators of $\fe_{C_{n}}$ act on the elements in the spanning set is given in Lemma \ref{typeccomputation}.\vspace{0.1in}

Let $S$ be the set $\{\eiaj{i}{j}|i<j\}\backslash\{[e_1[e_1e_2]]\}$. We have the following lemma:

\begin{lemma}
\label{commutativeideal}

Let $I'$ be the vector space spanned by $S$. Then $I'$ is an ideal of $\fe_{C_{n}}$.
\end{lemma}
\proof Based on Lemma \ref{typeccomputation}, we can see that $[e_{i},s]$ is a linear combination of elements in $S$ for all $s\in S$ and $i\in [n]$. $\square$\vspace{0.2in}

Furthermore, $I'$ has a special property:

\begin{lemma}
\label{Icommute}

The ideal $I'$ is abelian, that is, $[I',I']=0$.
\end{lemma}

\proof See Proof $\ref{proofIcommute}$.$\square$\vspace{0.2in}

We also need the following lemma for later.

\begin{lemma}
\label{prop1123}

\begin{alignat*}{3}
&[[e_1e_2],\one]=\one,                 \\
&[[e_3e_4],\one]=\one,                 \\
&[\eij{2i+1}{j-1}{j},\one]=0       &&\ \ \forall j >2i+1 \geq 3,j\neq 4,\\
&[\eij{1}{j-1}{j},\one]=0             &&\forall j\geq 3,\\
\end{alignat*}
\end{lemma}

\proof See Proof \ref{proofprop1123}.$\square$\vspace{0.2in}

Now consider the case when $n$ is even, that is, $\fe_{C_{2n}}$.

\begin{lemma}\label{centerofec}
Let 
\begin{align*}
c=&\ 2n\cdot e_1+n\cdot[e_1[e_1e_2]]\\&+\sum_{i=1}^{n-1}(n-i)(\eia{2i}{2i}+\eiai{2i})\\ 
&+\sum_{i=1}^{n-1}\sum_{j=1}^i(-1)^{i+j-1}\eiajthree{2j-1}{2i+1}{2i+2}. \\
\end{align*}
Then $c$ is in the center of $\fe_{C_{2n}}$.
\end{lemma}
\proof See Proof \ref{proofcenterofec}.$\square$\vspace{0.2in}

Define $I$ to be the vector space spanned by $S$ together with $c$. Lemma \ref{commutativeideal}, \ref{Icommute}, and \ref{centerofec} show that (1) $I$ is an ideal, and (2) $[I,I]=0$. Also, $\fe_{C_{2n}}/I$ is generated by $\bar{e}_i$'s via the relations $[\bar{e}_i[\bar{e}_i\bar{e}_{i\pm1}]]=-2\bar{e}_i$, for all $i$ except $i=1$, and $[\bar{e}_i,\bar{e}_j]=0$ for $|i-j|\ge 2$. However, the element $c$ in $I$ gives us the relation $[\bar{e}_1[\bar{e}_1\bar{e}_2]=-2\bar{e}_1$. This shows that $\fe_{c_{2n}}/I\cong \fe_{A_{2n}}\cong\spn_{2n}$.\vspace{0.1in}

Applying Lemma \ref{semisimpleaction} to our case, we see that $\fe_{C_{2n}}/I$ has an action on $I$. Our goal is to find how $I$ is decomposed into irreducible representations of $\spn_{2n}$, and show that $\fe_{C_{2n}}\cong \fe_{C_{2n}}/I\ltimes I\cong \spn_{2n}\ltimes I$.\vspace{0.1in}

To understand the structure of $\fe_{C_{2n}}$ we first have to understand the structure of the $\spn_{2n}$-representation $I$. The plan is to find the highest weight vectors in $I$.\vspace{0.1in}

First of all, let $E_{ij}$ be the $n\times n$ matrix whose $(i,j)$ entry is 1 and 0 otherwise. Thanks to \cite{lp}, we have an isomorphism $\phi$ from $\fe_{C_{2n}}/I\cong \fe_{A_{2n}}$ to $\mathfrak{sp}_{2n}$:
\begin{align*}
\bar{e}_{1}&\mapsto\begin{pmatrix}0& E_{11}\\0&0\end{pmatrix}, \quad\bar{e}_{2i-1}\mapsto \begin{pmatrix}0& E_{(i-1)(i-1)}+E_{(i-1)i}+E_{i(i-1)}+E_{ii}\\0&0\end{pmatrix} \text{ for $i\geq$2},\\
\bar{e}_{2i}&\mapsto \begin{pmatrix}0& 0\\E_{ii}&0\end{pmatrix}.
\end{align*}

We would like to find all of generators of $\fe_{C_{2n}}/I\cong \mathfrak{sp}_{2n}$ which correspond to simple roots, i.e. the preimage of $\begin{pmatrix}E_{(i-1)i}& 0\\0&-E_{i(i-1)}\end{pmatrix}$ for $i\geq 2$ and $\begin{pmatrix}0& E_{nn}\\0&0\end{pmatrix}$. Also we need to find the maximal toral subalgebra in $\fe_{C_{2n}}/I$, i.e. the preimage of $\begin{pmatrix}E_{ii}& 0\\0&-E_{ii}\end{pmatrix}$ for $i\geq 1$.

\begin{lemma}
\label{weightvectors}
\ 
\begin{enumerate}
\item \begin{align*}
          \phi^{-1}\begin{pmatrix}E_{kk}& 0\\0&-E_{kk}\end{pmatrix}=&\ \bbij{2k-1}{2k}-\beij{2k-3}{2k-1}{2k}\\
                                                                                                      &+\beij{2k-5}{2k-1}{2k}+\ldots+(-1)^{k-1}\beij{1}{2k-1}{2k},
          \end{align*}
\item \begin{align*}
          \phi^{-1}\begin{pmatrix}E_{(k-1)k}& 0\\0&-E_{k(k-1)}\end{pmatrix}=&\ \beij{2k-3}{2k-1}{2k}-\beij{2k-5}{2k-1}{2k}\\
                                                                                                      &+(-1)^k\beij{1}{2k-1}{2k},
         \end{align*}
\item \begin{align*}
 \phi^{-1}\begin{pmatrix}E_{(k+1)k}& 0\\0&-E_{k(k+1)}\end{pmatrix}=&\ \bbij{2k+1}{2k}-\bbij{2k-1}{2k}+\beij{2k-3}{2k-1}{2k}\\
                                                                                                      &-\beij{2k-5}{2k-1}{2k}+\ldots+(-1)^{k}\beij{1}{2k-1}{2k},
\end{align*}
\item \begin{align*}
         \phi^{-1}\begin{pmatrix}0& E_{nn}\\0&0\end{pmatrix}=&\ \bar{e}_1+\bar{e}_3+\ldots \bar{e}_{2n-1}\\  
                             &-[\bar{e}_1[\bar{e}_2\bar{e}_3]]-[\bar{e}_3[\bar{e}_4\bar{e}_5]]-\ldots-[\bar{e}_{2n-3}[\bar{e}_{2n-2}\bar{e}_{2n-1}]]\\
                             &+\beij{1}{4}{5}+\beij{3}{6}{7}+\ldots \beij{2n-5}{2n-2}{2n-1}\\
                             &+\ldots\\
                             &+(-1)^{n-1}\beij{1}{2n-2}{2n-1}\\
                          =&\sum_{l=0}^{n-1}\sum_{k=0}^{n-1-l}(-1)^l\beij{2k+1}{2k+2l}{2k+2l+1}.
\end{align*}
\end{enumerate}

\end{lemma}

\proof See Proof \ref{proofweightvectors}. $\square$\vspace{0.2in}

This computation leads to the following lemma:

\begin{lemma}
\label{annilator}

The elements $c$ and $\ee$ in $I$ are annihilated by $\phi^{-1}\begin{pmatrix}E_{(k-1)k}&0\\0&-E_{k(k-1)}\end{pmatrix}$ for all $k\geq 2$ and $\phi^{-1}\begin{pmatrix}0&E_{nn}\\0&0\end{pmatrix}$.
\end{lemma}

\proof

Since $c\in Z(\fe_{C_{2n}})$, the center of $\fe_{C_{2n}}$, clearly it is annihilated by the elements of $\fe_{C_{2n}}$. By Lemma \ref{prop1123}, we notice that when $k\geq 2$, the commutator $[\eij{2i-1}{2k-1}{2k},\ee]=0$ for all $i\leq k-1$. From Lemma \ref{weightvectors},we know that \[\phi^{-1}\begin{pmatrix}E_{(k-1)k}&0\\0&-E_{k(k-1)}\end{pmatrix}=\sum_{i=1}^{k-1}(-1)^{i+1}\beij{2(k-i)-1}{2k-1}{2k},\]
so the action of $\phi^{-1}\begin{pmatrix}E_{(k-1)k}&0\\0&-E_{k(k-1)}\end{pmatrix}$ on $\ee$ gives 0 for all $k\ge 2$. 

Similarly, by Lemma \ref{prop1123}, one knows $[\eij{2i-1}{2j-2}{2j-1},\ee]=0$ for all $i\leq j$. Again by Lemma \ref{weightvectors}, 
\[\phi^{-1}\begin{pmatrix}0& E_{nn}\\0&0\end{pmatrix}=\sum_{l=0}^{n-1}\sum_{k=0}^{n-1-l}(-1)^l\beij{2k+1}{2k+2l}{2k+2l+1},\]
 so the action of $\phi^{-1}\begin{pmatrix}0& E_{nn}\\0&0\end{pmatrix}$ on $\ee$ is 0. $\square$\vspace{0.1in}

By Lemma \ref{annilator}, it turns out that $c$ and $\ee$ are both highest weight vectors. We will find their weights. Because $c\in Z(\fe_{C_{2n}})$, its weight vector has to be a zero vector. Hence, the element $c$ spans a trivial representation of $\fe_{C_{2n}}/I\cong \mathfrak{sp}_{2n}$. As for $\ee$:\vspace{0.1in}

\begin{lemma}
\label{lambdaweight}

The weight of $\ee$ is $\omega_1+\omega_2$, where $\omega_1$ and $\omega_2$ are first and second fundamental weights of $\spn_{2n}$.
\end{lemma}

\proof

Apply Lemma \ref{prop1123}:
\[ [\bij{1}{2},\ee]=\ee,\]
so $\phi^{-1}\begin{pmatrix}E_{11}& 0\\0&-E_{11}\end{pmatrix}$ acts on $\ee$ by 1.

\begin{align*}
&[\bij{3}{4}-[e_1[e_2[e_3e_4]]],\ee]\\
=&[\bij{3}{4},\ee]-[[e_1[e_2[e_3e_4]]],\ee]\\
=&\ee-0=\ee,
\end{align*}
so $\phi^{-1}\begin{pmatrix}E_{22}& 0\\0&-E_{22}\end{pmatrix}$ acts on $\ee$ by 1.\vspace{0.1in}

For $k\geq 3$, $[\eij{2i-1}{2k-1}{2k},\ee]=0$ when $i\leq k$. Therefore, $\\ \phi^{-1}\begin{pmatrix}E_{kk}& 0\\0&-E_{kk}\end{pmatrix}=\sum_{i=0}^k(-1)^{i}\beij{2(k-i)-1}{2k-1}{2k}$ annihilates $\ee$. This completes the lemma. $\square$\vspace{0.2in}

Let $V_{\mathbf{0}}$ be the trivial $\spn_{2n}$-representation, and $V_{\lambda}$ be the irreducible $\spn_{2n}$-representation with highest weight $\lambda=\omega_1+\omega_2$. Lemma \ref{annilator} and \ref{lambdaweight} imply that $V_{\lambda}\oplus V_{\mathbf{0}}$ is isomorphic to an $\mathfrak{sp}_{2n}$-subrepresentation of $I$. By the Weyl character formula \cite{hj}, $\dim V_{\lambda}=(2n+1)(n-1)=2n^2-n-1$, and $\dim V_{\mathbf{0}}=1$, so $\dim V_{\lambda}+\dim V_{\mathbf{0}}=2n^2-n\le \dim I$. On the other hand, Since $S\cup \{c\}$ form a spanning set of $I$, $\dim I\le|S\cup\{c\}|=2n^2-n$. So $\dim I=\dim V_{\lambda}+\dim V_{\mathbf{0}}$. Thus $I\cong V_{\lambda}\oplus V_{\mathbf{0}}$.\vspace{0.2in}

The above argument is based on the assumption that $c$ and $\ee$ are not equal to 0. We still need to show $c$ and $\ee$ are not zero. \vspace{0.2in}

Let $F[i,j]$ be the $n^2\times n^2$ matrix with 1 in $i,j$ entry and 0 elsewhere. Define the a Lie algebra homomorphism (this is actually the adjoint representation of $\fe_{C_{2n}}$) from $\fe_{C_{2n}}$ to $\mathfrak{gl}_{n^2}$ by 
\begin{align*}
e_1\mapsto &F[3,2]+F[4,3]-F[8,16]+\sum_{j=3}^{2n}(F[(j-1)^2+j,(j-1)^2+j-1]\\
                    &F[(j-1)^2+j+1,(j-1)^2+j]+F[(j-1)^2+j+1,(j-1)^2+j+2]),
\end{align*}
\begin{align*}
e_k\mapsto &-\sum_{i=2}^{2k-2}F[(k-1)^2+i,(k-2)^2+i-1]+2F[(k-1)^2+1,(k-1)^2+2]\\
                    &+\begin{cases}F[(k-1)^2+2k-2,(k-1)^2+2k-1]& \text { if }k\ge 3\\2F[3,4]&\text{ if } k=2\end{cases}\\
                    &+F[k^2+2,k^2+1]-2F[k^2+2,k^2+3]-\sum_{i=3}^{2k-1}F[(k-1)^2+i-1,k^2+i]\\
                    &+F[k^2+2k+1,k^2+2k]-F[(k+1)^2+2k+2,(k+2)^2+2k+5]\\
                    &+\sum_{j=k+2}^{2n}(F[(j-1)^2+j-k+1,(j-1)^2+j-k]\\
                    &\ \ \ \ \ \ \ +F[(j-1)^2+j-k+1,(j-1)^2+j-k+2]\\
                    &\ \ \ \ \ \ \  +F[(j-1)^2+j+k,(j-1)^2+j+k-1]\\
                    &\ \ \ \ \ \ \  +F[(j-1)^2+j+k,(j-1)^2+j+k+1]).
\end{align*}

It is straightforward to verify this map satisfies all of the relation among $e_k$, thus a Lie algebra homomorphism of $\fe_{C_{2n}}$. And the $(14,10)$ entry of the image of $[e_1[e_1[e_2e_3]]]$ is 1. Hence, $[e_1[e_1[e_2e_3]]]$ is not zero.\vspace{0.1in}

On the other hand, if we consider the map from $\fe_{C_{2n}}$ to $\mathfrak{gl}_1$:
\[e_1\mapsto 1,\qquad e_k\mapsto 0.\]

It will also satisfy the relation among $e_k's$, and the image of $c$ is $2n$, so	 $c$ is not zero.\vspace{0.1in}

\begin{theorem}
\label{typecstructure}
We have
\[\fe_{C_{2n}}\cong \spn_{2n}\ltimes(V_{\lambda}\oplus V_{\mathbf{0}}),\] 
where $\lambda=(1,1,0\ldots,0)$, $\mathbf{0}=(0,0,0\ldots,0)$.
\end{theorem}

\proof






We will use the second Lie algebra cohomology group $H^2(L, V)$, where $L$ is a Lie algebra, and $V$ is a representation of $L$. It is known that $H^2(L,V)$ is in bijection with extensions of $L^*$ with abelian kernal $V$ \cite{ce}. In our case, $\fe_{C_{2n}}$ is an extension of $\fe_{C_{2n}}/I\cong\fe_{A_{2n}}$ with abelian kernal $V\cong V_{\lambda}\oplus V_{\mathbf{0}}$. By Theorem 26.3 and 26.4 of \cite{ce}, since $\fe_{A_{2n}}\cong\spn_{2n}$ is semisimple, and $V$ is a finite dimensional representation of $\spn_{2n}$, we know that $H^2(\spn_{2n},V)=0$. So there is only one extension up to isomorphism, that is, $\spn_{2n}\ltimes V$. Hence, \[\fe_{C_{2n}}\cong \spn_{2n}\times (V_{\lambda}\oplus V_{\mathbf{0}}),\]
where $V_{\lambda}$ is the irreducible representation of $\spn_{2n}$ with the highest weight vector $\lambda=\omega_1+\omega_2$, and $V_{\mathbf{0}}$ is the trivial representation .\vspace{0.1in}

 One immediate corollary is:

\begin{corollary}
The dimension of $\fe_{C_{n}}$ is $n^2$.
\end{corollary}
\proof $\dim \fe_{C_{2n}}=\dim \fe_{A_{2n}}+\dim V_{\lambda}+\dim V_{\mathbf{0}}=2n^2+n+2n^2-n=4n^2$. The spanning set of $\fe_{C_{2n+1}}$ of our choice is a subset of a basis of $\fe_{C_{2n+2}}$, so they have to be linearly independent, thus a basis of $\fe_{C_{2n+1}}$, so $\dim \fe_{C_{2n+1}}=(2n+1)^2$.$\square$\vspace{0.2in}

\subsection{Type $D$}
\label{typed}

We will study the case of odd rank first. Electrical Lie algebra $\fe_{D_{2n+1}}$ is generated by generators $\{\dop,\don,e_2,e_3,\ldots,e_{2n}\}$ with the relations (Figure \ref{typeddynkin}):

\begin{figure}
\begin{tikzpicture}[scale=1.4]
\fill

(1,0.5) node(1){} circle(2pt)
(1,-0.5) node(1b){} circle(2pt)
(2,0) node(2){}circle(2pt)
(3,0) node(3){} circle(2pt)
(4,0) node(4){} circle(2pt)
(5,0) node(5){} circle(2pt)
(6,0) node(6){} circle(2pt)
(7,0) node(7){} circle(2pt)
(1,-0.5) node(n1)[below]{1}
(1,0.5) node(n1b)[above]{$\overline{1}$}
(2,0) node(n2)[above]{2}
(3,0) node(n3)[above]{3}
(4,0) node(n4)[above]{4}
(5,0) node(n5)[above]{$n-2$}
(6,0) node(n6)[above]{$n-1$}
(7,0) node(n7)[above]{$n$}
;

\draw (1) -- (2);
\draw (1b)--(2);
\draw (2) -- (3);
\draw (3) -- (4);
\draw[dashed] (4) -- (5);
\draw (5) -- (6);
\draw (6) -- (7);

\end{tikzpicture}
\caption{Dynkin Diagram of $\cD_{n+1}$}
\label{typeddynkin}
\end{figure}

\begin{alignat*}{3}
&[\dop,[\dop,e_2]]=-2\dop, \\
&[\don,[\don,e_2]]=-2\don,\\
&[e_2,[e_2,\dop]]=-2e_2, \\
&[e_2,[e_2,\don]]=-2e_2,\\
&[\dop,e_{i}]=  0 &&\text{ if } i\ge 3,\\
&[\don,e_{i}]=0  &&\text{ if } i\ge 3,\\
&[e_i,[e_i,e_j]]=-2e_i &&\text{ if } |i-j|=1, \\
&[e_i,e_j]=        0 &&\text{ if } |i-j|\ge 2,\ \ i,j\ge 2.
\end{alignat*}

By Proposition \ref{dimupperbound}, $\fe_{D_{2n+1}}$ has a spanning set:
\begin{align*}
&\{\dop,\don\}\cup\{\deij{i}{j-1}{j}\}_{2\le i\le j\le 2n}\cup\{\dei{j-1}{j}\}_{2\le j\le 2n}\\
&\cup\{\depi{j-1}{j}\}_{2\le j\le 2n}\cup\{\dejj{j-1}{j}\}_{2\le j\le 2n}\\
&\cup\{\deiajthree{i}{j-1}{j}\}_{2\le i<j\le 2n}.
\end{align*}
The brackets of generators with the elements in the spanning set are entirely similar to type $C$, which we will omit here.\vspace{0.1in}

Let 
\begin{align*}
c
=&\ n\cdot e_1+n\cdot \dop+n\cdot[\dop[e_1e_2]]\\&+\sum_{i=1}^{n-1}(n-i)(\deiajthree{2i}{2i}{2i+1}+\deiajthree{2i+1}{2i+1}{2i+2})\\ 
&+ \sum_{i=1}^{n-1}\sum_{j=1}^i(-1)^{i+j-1}\deiajthree{2j-1}{2i+1}{2i+2}
\end{align*}\vspace{0.1in}

\begin{lemma}
\label{ecnsubalgebra}
Let $\{f_i\}_{i\ge 1}$ be the generating set of $\fe_{C_{2n}}$. Consider the map $\phi: \fe_{C_{2n}}\rightarrow \fe_{D_{2n+1}}$:
\[\phi(f_1)=\frac{e_1+\dop}{2},\qquad \phi(f_k)=e_k\ \  \forall k\ge 2\]
 and extend it by linearity. Then $\phi$ is a Lie algebra homomorphism. 
\end{lemma}

\proof We only need to show the relations $[\frac{e_1+\dop}{2},[\frac{e_1+\dop}{2},[\frac{e_1+\dop}{2},e_2]]]=0$, $[e_2,[e_2,\frac{e_1+\dop}{2}]]=-2e_2$, and $[\frac{e_1+\dop}{2},e_k]=0$ for $k\ge 3$ in $\fe_{D_{2n+1}}$. These can all be verified by straightforward computation.$\square$ \vspace{0.2in}

Recall that the radical ideal $I'$ of $\fe_{C_{2n}}$ is an abelian ideal which is a direct sum of  its center 
\begin{align*}
c'=&\ 2n\cdot f_1+n\cdot[f_1[f_1f_2]]\\
     &+\sum_{i=1}^{n-1}(n-i)(\deiajthreec{2i}{2i}{2i+1}+\deiajthreec{2i+1}{2i+1}{2i+2})\\
     &+ \sum_{i=1}^{n-1}\sum_{j=1}^i(-1)^{i+j-1}\deiajthreec{2j-1}{2i+1}{2i+2}.
\end{align*}
 and an abelian ideal $I_1'$ with basis vector of the form $\deiajthreec{i}{j-1}{j}$, where $i<j$, $j\ge 3$. Now we calculate the image of $I'$ under the map $\phi$. We will use an identity for $j\ge 3$, $\deiajthreet{i}{j-1}{j}=\frac{1}{2}\deiajthree{i}{j-1}{j}$, and $[\frac{e_1+\dop}{2},[\frac{e_1+\dop}{2},e_2]]=-\frac{1}{2}(e_1+\dop)+\frac{1}{2}[\dop[e_1e_2]]$.\vspace{0.1in}

Thus,
\begin{align*}
&\phi(\deiajthreec{i}{j-1}{j})=\frac{1}{2}\deiajthree{i}{j-1}{j}\quad \forall j\ge 3,\\
&\phi(c')= n(e_1+\dop)-\frac{n}{2}(e_1+\dop)+\frac{n[\dop[e_1e_2]]}{2}\\
&\quad\qquad +\frac{1}{2}(\sum_{i=1}^{n-1}(n-i)(\deiajthree{2i}{2i}{2i+1}+\deiajthree{2i+1}{2i+1}{2i+2})\\ 
&\quad\qquad + \sum_{i=1}^{n-1}\sum_{j=1}^i(-1)^{i+j-1}\deiajthree{2j-1}{2i+1}{2i+2})\\
&\qquad=\frac{1}{2}c.
\end{align*}

\begin{lemma}
\label{Iimage}
Let $I$ be the image $\phi(I')$. Then $I$ is an ideal of $\fe_{D_{2n+1}}$. Moreover, $\phi$ is injective.
\end{lemma}

\proof Let $y=\phi(x)\in I$. Clearly, for $k\ge 2$, $[e_k,y]=\phi([f_k,x])\in I$ and $[e_k,c]=\phi([f_k,c'])=0\in I$. It suffices to show that $[\dop,y]\in I$ and $[e_1,y]\in I$.\vspace{0.1in}

If $y=\deiajthree{i}{j-1}{j}$ where $j\ge 3,\  i<j$. By straightforward computation, we have 

\begin{align} \label{e1action}
[e_1,y]=\begin{cases}\dejj{j-1}{j} &\text{ if } i=2,\\ 
                                        [e_{\bar{1}}[e_1[e_2e_3]]] &\text{ if }i=3, j=4,\\
                                         0 &\text{ otherwise. }\end{cases} 
\end{align}
By symmetry, we will get the same results as above for $\dop$. So $[e_1,y]=[e_{\bar{1}},y]\in I$.\vspace{0.1in}

If $y=c$, we know that $[(e_1+e_{\bar{1}})/2,c]=\phi([f_1,c'])=0$. From Equation \eqref{e1action}, $[e_1-e_{\bar{1}},c]=0$. Thus $[e_1,c]=[e_{\bar{1}},c]=0\in I$.\vspace{0.1in}

Therefore, $I$ is an ideal of $\fe_{D_{2n+1}}$. Furthermore, the above shows that $c$ is in the center of $\fe_{D_{2n+1}}$.\vspace{0.1in}

Consider the image $\phi(\fe_{C_{2n}})$. It naturally has a spanning set which is the image of basis elements of $\fe_{C_{2n}}$ under $\phi$. If we compute the Lie bracket of this spanning set with itself, it is the same as the adjoint representation of $\fe_{C_{2n}}$, thus has dimension $4n^2-1$. It suffices to show that $c$ is not zero. Let $I_1=\phi(I_1')$, which is spanned by the elements of the form $\deiajthree{i}{j-1}{j}$, where $j\ge 3$. By above computation, it is an ideal of $\fe_{D_{2n+1}}$. In $\fe_{D_{2n+1}}/I_1$, we have $\bar{c}=n(\bar{e}_1+\bar{e}_{\bar{1}}+[\bar{e}_{\bar{1}}[\bar{e}_1\bar{e}_2])$. Note that $\bar{e}_1$, $\bar{e}_2$, and $\bar{e}_{\bar{1}}$ form a Lie subalgebra of $\fe_{D_{2n+1}}/I$, which is isomorphic to $\fe_{A_3}$. By the isomorphism in the proof of Theorem \ref{oddtypeastructure1} $\bar{e}_1+\bar{e}_{\bar{1}}+[\bar{e}_{\bar{1}}[\bar{e}_1\bar{e}_2]$ is not zero, and in the center of $\fe_{A_3}$. Thus $\bar{c}$ is not zero, neither is $c\in \fe_{D_{2n+1}}$. So we have $\dim \phi(\fe_{C_{2n}})=4n^2=\dim \fe_{C_{2n}}$. We conclude that $\phi$ is injective.$\square$\vspace{0,2in}

Let $J$ be the ideal generated by $\dop-\don$. It is clear that $\fe_{D_{2n+1}}/J$ is isomorphic to $\fe_{A_{2n}}$. We study the structure of the ideal $J$. Again by Equation \eqref{e1action} and computation of type $A$ electrical Lie algebra, a spanning set for $J$ is 
\[\{e_1-\dop,c\}\cup\{[e_1-\dop[e_2[\ldots[e_{j-1}e_j]\ldots]\}_{j=2}^{2n}\cup \{\deiajthreetm{i}{j-1}{j}\}_{j\ge 2,i<j}.\] 
Note that $\deiajthreetm{i}{j-1}{j}=-2\deiajthree{i}{j-1}{j}$ for $j\ge 3$. The reason why $c\in J$ is because 
\[[e_1-\dop,[e_1-\dop,e_2]]=-2e_1-2\dop-2[\dop[e_1e_2]]=\frac{2}{n}(c-\text{ linear combination of } \deiajthree{i}{j-1}{j}),\] where $j\ge 3$. By this observation, we have $I\subsetneq J$. Let $K=\{e_1-\dop,[e_1-\dop,e_2],\ldots,[e_1-\dop,[e_2[\ldots[e_{2n-1}e_{2n}]\ldots]\}$. Then $J$ is spanned by $I$ and $K$ as a vector space.\vspace{0.1in}

We claim that $[J,J]\in I$. If this is true, then because $[I,I]=0$ from Lemma \ref{Icommute} and \ref{Iimage}, $J$ is the radical of $\fe_{D_{2n+1}}$.\vspace{0.1in}

To show the above claim, it suffices to find $[K,I]$ and $[K,K]$. Due to Equation $\eqref{e1action}$, $[e_1-\dop,I]=0$. Hence by simple induction on $k$, we have $[[e_1-\dop,[e_2[\ldots[e_{k-1}e_k]\ldots],I]=0$. Hence $[K,I]=0$ 

\begin{lemma} 
\label{typedcomputation}

Assume $i\ge j$. Then
\begin{align*}
&[\dei{i-1}{i},\dei{j-1}{j}]=\begin{cases}2e_1+[e_1[e_2e_3]]& \text{ if } i=2, j=3,\\
                           (-1)^i2e_1                 & \text{ if } j=i+1,i\neq 2,\\
                            0                               &\text{ if }  |j-i|\ge 2 \text{ or } j=i,\end{cases}\\
&[\depi{i-1}{i},\depi{j-1}{j}]=\begin{cases}2\dop+[\dop[e_2e_3]]& \text{ if } i=2, j=3,\\
                           (-1)^i2\dop                 & \text{ if } j=i+1,i\neq 2,\\
                           0                                &\text{ if }  |j-i|\ge 2 \text{ or } j=i,\end{cases}\\
&[\depi{i-1}{i},\dei{j-1}{j}]\\
=&\begin{cases}-[e_2[\dop[e_1[e_2e_3]]]]-[\dop[e_1e_2]]+[\dop[e_2e_3]]& \text{ if } i=2, j=3,\\
                            (-1)^{i-1}([\dop[e_1e_2]]+\sum_{s=2}^i\deiajthree{s}{s}{s+1})      & \text{ if } j=i+1, j\ge 4,\\
                            (-1)^{i-1}\deiajthree{i}{j-1}{j}                                 &\text{ if }  |j-i|\ge 2,\\
                            [e_{\bar{1}}e_2]-[e_1e_2]              &\text{ if } i=j=2,\\
                            0                                                       &\text{ if } i=j\ge 3.\end{cases}\\
\end{align*}
\end{lemma}

\proof See Proof \ref{proofsdcomputation}.$\square$\vspace{0.2in}

By Lemma \ref{typedcomputation}, we obtain that 
\begin{align*}
&[[e_1-\dop,[e_2[\ldots[e_{i-1}e_i]\ldots],[e_1-\dop,[e_2[\ldots[e_{j-1}e_j]\ldots]]\\
=&\begin{cases} (-1)^i(2(e_1+\dop)+2[\dop[e_1e_2]]+2\sum_{s=2}^{i}\deiajthree{s}{s}{s+1}) &\text{ if } j=i+1,\\
                             (-1)^i2\deiajthree{i}{j-1}{j}                          &\text{ if } |j-i|\ge 2.
\end{cases}
\end{align*}

Since $2(e_1+\dop)+2[\dop[e_1e_2]]$ is a linear combination of $c$ and $\deiajthree{i}{j-1}{j}$ where $j\ge 3,\  i<j$, we have $[K,K]\in I$. Combining the above calculation, we prove $[J,J]\in I$.

Then we have the folllowing theorem:\vspace{0.1in}

\begin{theorem}
$J$ is the radical ideal of $\fe_{D_{2n+1}}$. Furthermore, $\fe_{D_{2n+1}}/J\cong \spn_{2n}$.
\end{theorem}\vspace{0.1in}

As a consequence of the above calculation, we have:\vspace{0.1in}

\begin{theorem}
\label{typeddimension} We have
\[\dim \fe_{D_{n}}=n^2-n.\]

\end{theorem}

\proof First consider the case when $n$ is odd. Using the above notations, we know $[K,K]\in I$, so for $\bar{K}$ as an ideal in $\fe_{D_{2n+1}}/I$, we have $[\bar{K},\bar{K}]=\bar{0}$. Thus, by Lemma \ref{semisimpleaction}, $\bar{K}$ is a representation of the quotient Lie algebra $(\fe_{D_{2n+1}}/I)/\bar{K}$. Note that $(\fe_{D_{2n+1}}/I)/\bar{K}\cong \fe_{D_{2n+1}}/J\cong \spn_{2n}$. \vspace{0.1in}

We claim that $\bar{K}\neq \bar{0}$. Otherwise, the injective homomorphism $\phi$ defined in Lemma \ref{ecnsubalgebra} is an isomorphism. Then the  element $e_1$ is not in $\fe_{D_{2n+1}}$, a contradiction. \vspace{0.1in}

Furthermore, $\bar{K}$ is cyclic generated by $\bar{e}_1-\bar{e}_{\bar{1}}$. Thus $\bar{K}$ is irreducible. Because $\dim \bar{K}\leq 2n$ and the nontrivial irreducible representation of $\spn_{2n}$ with the smallest dimension is the standard representation $V_{\nu}$, where $\nu=\omega_1$, the first fundamental weight, we have that $\bar{K}\cong V_{\nu}$. On the other hand, by Lemma \ref{Iimage} the homomorphism $\phi: \fe_{C_{2n}}\longrightarrow \fe_{D_{2n+1}}$ is injective, so the spanning set of $I$ is indeed a basis. Hence
\[\dim \fe_{D_{2n+1}}=\dim \spn_{2n}+\dim V_{\nu}+\dim I=2n^2+n+2n+2n^2-n=(2n+1)^2-(2n+1)\]

Since $\fe_{D_{2n}}$ is a Lie subalgebra of $\fe_{D_{2n+1}}$, by an argument similar to the one for type $C$, we can prove that $\dim \fe_{D_{2n}}=(2n)^2-(2n)$. $\square$ \vspace{0.1in}

We also have the following conjecture:

\begin{conjecture}
\label{typedstructure}

We conjecture that $\fe_{D_{2n+1}}$ is isomorphic to an extension of $\spn_{2n}\ltimes V_{\nu}$ by $V_{\lambda}\oplus V_{\mathbf{0}}$, where $V_{\nu}$ is the standard representation, $V_{\mathbf{0}}$ is the trivial representation, and $V_{\lambda}$ is the irreducible representation of $\spn_{2n}$ with highest weight vector $\lambda$ being the sum of the first and second fundamental weights. In another word, there is a short exact sequence
\[0\longrightarrow V_{\lambda}\oplus V_{\mathbf{0}}\longrightarrow \fe_{D_{2n+1}}\longrightarrow \spn_{2n}\ltimes V_{\nu}\longrightarrow 0.\]
\end{conjecture}\vspace{0.2in}

\section{Facts and Proofs of Lemmas}

\subsection{Facts and Proofs for Type $C$ in \ref{typec}}
\label{typecproof}

The following lemma describes how the generators $e_i$ act on the spanning set of $\fe_{C_{n}}$.

\begin{lemma}\label{typeccomputation}

Bracket $e_k$ with the elements in the spanning set:\vspace{0.1in}

$[e_k, \eij{i}{j-1}{j}]$:\vspace{0.05in}

If $j=i$ and $j=i+1$ it is given in the defining relation; if $k<i-1$ or $k>j+1$, $[e_k,\eij{i}{j-1}{j}]=0$; if $j\geq i+2$ and $i-1\leq k\leq j+1$:

\begin{enumerate}
\item[(1)] $j=i+2$
\begin{align*}
&[e_{i-1},\eijk{i}{i+1}{i+2}]=\eij{i-1}{i+1}{i+2},\\
&[e_i,\eijk{i}{i+1}{i+2}]=\begin{cases}0 &\text{ if } i\neq 1,\\ \eijkr{1}{1}{2}{3}
&\text{ if } i=1,\end{cases}\\
&[e_{i+1},\eijk{i}{i+1}{i+2}]=-\bij{i}{i+1}+\bij{i+1}{i+2},\\
&[e_{i+2},\eijk{i}{i+1}{i+2}]=0,\\
&[e_{i+3},\eijk{i}{i+1}{i+2}]=-\eijkr{i}{i+1}{i+2}{i+3}.
\end{align*}

\item[(2)] $j\geq i+3$

\begin{align*}
&[e_{i-1},\eij{i}{j-1}{j}]=\eij{i-1}{j-1}{j},\\
&[e_{i},\eij{i}{j-1}{j}]=\begin{cases} 0 &\text{  if }i\neq 1, \\ [e_1[e_1[\ldots[e_{j-1}e_j]\ldots] &\text{  if } i=1,\end{cases}\\
&[e_{i+1},\eij{i}{j-1}{j}]=\eij{i+1}{j-1}{j},\\
&[e_k,\eij{i}{j-1}{j}]=0 \text{ if $i+2\leq k\leq j-2$ },\\
&[e_{j-1},\eij{i}{j-1}{j}]=-\eij{i}{j-2}{j-1},\\
&[e_j,\eij{i}{j-1}{j}]=0,\\
&[e_{j+1},\eij{i}{j-1}{j}=-\eij{i}{j}{j+1}.
\end{align*}
\end{enumerate}\vspace{0.2in}

$[e_k,\eiajj{i}{i-1}{j-1}{j}]$, where $j\geq i+1$:\vspace{0.05in}

If $k\ge j+2$, $[e_k,\eiajj{i}{i-1}{j-1}{j}=0$. If $k\leq {j+1}$:

\begin{enumerate}
\item[(3)] $j\geq i+2$
\begin{align*}
&[e_{j+1},\eiajj{i}{i-1}{j-1}{j}]=-\eiaj{i}{j},\\
&[e_j,\eiaj{i}{j}]=0,\\
&[e_{j-1},\eiaj{i}{j}]=-\eiajfive{i}{i-1}{j-3}{j-2}{j-1},\\
&[e_k,\eiaj{i}{j}]=0 \text{ if $i+2\leq k\leq j-2$},\\
&[e_{i+1},\eiaj{i}{j}]= \eiajj{i+1}{i}{j-1}{j},\\
&[e_i,\eiaj{i}{j}]=0,\\
&[e_{i-1},\eiajj{i}{i-1}{j-1}{j}]=\eiajj{i-1}{i-2}{j-1}{j},\\
&[e_k,\eiaj{i}{j}]=0\text{ if $k\leq i-2$}.
\end{align*}

\item[(4)]  $j=i+1$

\begin{align*}
&[e_{i+2},\eiajthree{i}{i}{i+1}]=-\eiajthree{i}{i+1}{i+2},\\
&[e_{i+1},\eia{i}{i}]=\begin{cases} \eia{i-1}{i} &\text{ if } i \geq 2, \\  2[e_1e_2] &\text{ if } i=1,\end{cases}\\
&[e_{i},\eia{i}{i}]=0,\\
&[e_{i-1},\eiajj{i}{i-1}{i}{i+1}]=\eia{i-1}{i},\\
&[e_{i-2},\eia{i}{i}]=-\eiaj{i-2}{i},\\
&[e_k,\eia{i}{i}]=0 \text{ if $ k\le i-3$}.
\end{align*}
\end{enumerate}

\end{lemma}
\proof 
The above identities can be achieved by the straightforward but lengthy computation. We will omit it here.\vspace{0.2in}

\begin{proofs}[of Lemma \ref{Icommute}]
\label{proofIcommute}

\end{proofs}
We would like to show $[\eiaj{i_1}{j_1},\eiaj{i_2}{j_2}]=0$, where $i_1<j_1,i_2<j_2$ ,and $j_1,j_2\geq 3$. This will be done by induction on $i_1+j_1+i_2+j_2$. Let $i_k+j_k$ be the \textit{length} of $\eiaj{i_k}{j_k}$ in $S$. First, we show a few base cases which will be used shortly.

Base cases:

\begin{enumerate}

\item \[[[e_1[e_1[e_2[e_3e_4]]]],[e_1[e_1[e_2e_3]]]]=0.\]
First of all
\begin{align*}
&[[e_2e_3],[e_1[e_1[e_2[e_3e_4]]]]]\\ 
=&\ [[e_2[e_1[e_1[e_2[e_3e_4]\ldots],e_3]+[e_2[e_3[e_1[e_1[e_2[e_3e_4]\ldots] \\
=&-[e_3[e_2[e_1[e_1[e_2[e_3e_4]\ldots]-[e_2[e_1[e_1[e_2e_3]]]],
\end{align*}
so
\begin{align*}
&[[e_1[e_1[e_2[e_3e_4]]]],[e_1[e_2e_3]]]\\
=&\ [[e_2e_3],[e_1[e_1[e_1[e_2[e_3e_4]\ldots]]-[e_1,[[e_2e_3][e_1[e_1[e_2[e_3e_4]\ldots]]]\\
=&\ [e_1[e_3[e_2[e_1[e_1[e_2[e_3e_4]\ldots]+[e_1[e_2[e_1[e_1[e_2e_3]\ldots]\\
=&-[e_1[e_1[e_2e_3]]]+[e_1[e_1[e_2e_3]]]\\
=&\ 0.
\end{align*}
Thus,
\begin{align*}
&[[e_1[e_1[e_2[e_3e_4]]]],[e_1[e_1[e_2e_3]]]]\\
=&\ [[e_1[e_2e_3]],[e_1[e_1[e_1[e_2[e_3e_4]\ldots]]-[e_1,[ [e_1[e_2e_3]][e_1[e_1[e_2[e_3e_4]\ldots]]]\\
=&\ 0.\\
\end{align*}

\item \[[[e_1[e_1[e_2e_3]]],[e_2[e_1[e_1[e_2e_3]]]]]=0.\]

First we compute
\begin{align*}
&[[e_2e_3],[e_1[e_1[e_2e_3]]]]\\
=&-[e_3[e_2[e_1[e_1[e_2e_3]\ldots]+[e_2[e_3[e_1[e_1[e_2e_3]\ldots]\\
=&-[e_1[e_1[e_2e_3]]],
\end{align*}
and
\begin{align*}
&[[e_1[e_2e_3]],[e_1[e_1[e_2e_3]]]]\\
=&-[[e_3[e_1e_2]],[e_1[e_1[e_2e_3]]]]\\
=&\ [[e_1e_2],[e_3[e_1[e_1[e_2e_3]]]]]-[e_3,[[e_1e_2][e_1[e_1[e_2e_3]]]]]\\
=&\ 0-[e_3[e_1[e_1[e_2e_3]]]]\\
=&\ 0.
\end{align*}
Then,
\begin{align*}
&[[e_2[e_1[e_1[e_2e_3]]]],[e_1[e_2e_3]]]\\
=&\ [[e_2[e_1[e_2e_3]]],[e_1[e_1[e_2e_3]]]]+[[[e_1[e_2e_3]][e_1[e_1[e_2e_3]]]],e_2]\\
=&-[[e_1e_2][e_1[e_1[e_2e_3]]]]+[[e_2e_3],[e_1[e_1[e_2e_3]]]]+0\\
=&-2[e_1[e_1[e_2e_3]]].
\end{align*}
Therefore,
\begin{align*}
&[[e_1[e_1[e_2e_3]]],[e_2[e_1[e_1[e_2e_3]]]]]\\
=&\ [[e_1[e_2[e_1[e_1[e_2e_3]\ldots],[e_1[e_2e_3]]]+[[[e_2[e_1[e_1[e_2e_3]]]][e_1[e_2e_3]]],e_1]\\
=&\ [[e_1[e_1[e_2e_3]]],[e_1[e_2e_3]]]+2[e_1[e_1[e_1[e_2e_3]]]]\\
=&\ 0.
\end{align*}

\item \[[[e_1[e_1[e_2e_3]]],[e_2[e_1[e_1[e_2[e_3e_4]\ldots]]=0.\]
First we compute
\begin{align*}
&[[e_2e_3],[e_2[e_1[e_1[e_2[e_3e_4]\ldots]]\\
=&\ [[e_2[e_2[e_1[e_1[e_2[e_3e_4]\ldots],e_3]+[e_2[e_3[e_2[e_1[e_1[e_2[e_3e_4]\ldots]\\
=&\ 0+[e_2[e_1[e_1[e_2[e_3e_4]\ldots]\\
=&\ [e_2[e_1[e_1[e_2[e_3e_4]\ldots].
\end{align*}
Then
\begin{align*}
&[[e_1[e_2e_3]],[e_2[e_1[e_1[e_2[e_3e_4]\ldots]]\\
=&\ [e_1,[[e_2e_3][e_2[e_1[e_1[e_2[e_3e_4]\ldots]]]-[[e_2e_3],[e_1[e_2[e_1[e_1[e_2[e_3e_4]\ldots]]\\
=&\ [e_1[e_2[e_1[e_1[e_2[e_3e_4]\ldots]-\underbrace{[[e_2e_3],[e_1[e_1[e_2[e_3e_4]]]]]}_{\text{by an equation in case (1)}}\\
=&\ [e_1[e_1[e_2[e_3e_4]]]]+[e_3[e_2[e_1[e_1[e_2[e_3e_4]\ldots]+[e_2[e_1[e_1[e_2e_3]\ldots].\\
\end{align*}
Hence,
\begin{align*}
&[[e_1[e_1[e_2e_3]]],[e_2[e_1[e_1[e_2[e_3e_4]\ldots]]\\
=&\ [[e_1[e_2[e_1[e_1[e_2[e_3e_4]\ldots],[e_1[e_2e_3]]]\\
  &+[e_1,[[e_1[e_2e_3]][e_2[e_1[e_1[e_2[e_3e_4]\ldots]]]\\
=&\ \underbrace{[[e_1[e_1[e_2[e_3e_4]]]],[e_1[e_2e_3]]]}_{\text{by an equation in case (1)}}+[e_1[e_1[e_1[e_2[e_3e_4]\ldots]\\
&+[e_1[e_3[e_2[e_1[e_1[e_2[e_3e_4]\ldots]+[e_1[e_1[e_1[e_2e_3]\ldots]\\
=&\ 0+0-[e_1[e_1[e_2e_3]]]+[e_1[e_1[e_2e_3]]]\\
=&\ 0.
\end{align*}
\end{enumerate}

This finishes the base case we will use. Next we proceed to induction step. Up to symmetry, there are three cases: (i) $i_1\leq i_2<j_2\leq j_1$, (ii) $i_1\leq i_2\leq j_1<j_2$, where two equalities cannot achieve at the same time, (iii) $i_1<j_1\leq i_2<j_2$. We will use underbrace to indicate where induction hypothesis is used.

\begin{enumerate}

\item [(i)] If $i_1\leq i_2<j_2\leq j_1$
\begin{enumerate}
\item[\textbullet] If $i_1>1$
\begin{align*}
&[\eiaj{i_1}{j_1},\eiaj{i_2}{j_2}]\\
=&\ [\underbrace{[e_{i_1}\eiaj{i_2}{j_2}}_{\text{length decreases}},\eiaj{i_1-1}{j_1}]\\
&+[\underbrace{[\eiaj{i_2}{j_2},\eiaj{i_1-1}{j_1}]},e_{i_1}]\\
=&\ 0.
\end{align*}

\item [\textbullet] If $i_1=1$,and $j_1-j_2\geq 2$
\begin{align*}
&[\eiaj{i_1}{j_1},\eiaj{i_2}{j_2}]\\
=&-[[e_{j_1}\eiajthree{i_1}{j_1-2}{j_1-1}],\eiaj{i_2}{j_2}\\
=&-[\underbrace{[e_{j_1}\eiaj{i_2}{j_2}]}_{\text{equals } 0},\eiajthree{i_1}{j_1-2}{j_1-1}]\\
&-[\underbrace{[\eiaj{i_2}{j_2}\eiaj{i_1}{j_1-1}]}_{\text{length decreases and $j_1-1\geq j_2+1\geq 3$}},e_{j_1}]\\
=&\ 0.
\end{align*}

\item [\textbullet] If $i_1=1$, $j_1-j_2=1$, $j_2-i_2\geq 2$, and $j_2\geq 4$
\begin{align*}
&[\eiaj{i_1}{j_1},\eiaj{i_2}{j_2}]\\
=&-[\eiaj{i_1}{j_1},[e_{j_2}\eiajthree{i_2}{j_2-2}{j_2-1}]\\
=&-[\eiajthree{i_2}{j_2-2}{j_2-1},\underbrace{[e_{j_2}\eiaj{i_1}{j_1}}_{\text{length decreases}}]\\
&-[e_{j_2},\underbrace{[\eiaj{i_1}{j_1}\eiajthree{i_2}{j_2-2}{j_2-1}]}]\\
=&\ 0.
\end{align*}

\item [\textbullet] If $i_1=1$, $j_1-j_2=1$, $j_2-i_2\geq 2$, and $j_2= 3$, then we have to have $i_1=1,i_2=1,j_1=4,j_2=3$, this is base case (1).\vspace{0.1in}

\item [\textbullet] If $i_1=1$, $j_1-j_2=1$, $j_2-i_2=1$, and $i_2\geq 3$, then
\begin{align*}
&[\eiaj{i_1}{j_1},\eiaj{i_2}{j_2}]\\
=&\ [\underbrace{\eiaj{i_2-1}{j_2}}_{\text{length decreases}},\underbrace{[e_{i_2}\eiaj{i_1}{j_1}}_{\text{length decreases}}]\\
&+[e_{i_2},\underbrace{[\eiaj{i_1}{j_1}\eiaj{i_2-1}{j_2}]}]\\
=&\ 0.
\end{align*}

\item [\textbullet] If $i_1=1$, $j_1-j_2=1$, $j_2-i_2=1$, and $i_2=2$, then $i_1=1,i_2=2,j_1=4,j_2=3$
\begin{align*}
&[e_1[e_1[e_2[e_3e_4]\ldots],[e_2[e_1[e_1[e_2e_3]\ldots]]\\
=&\ \underbrace{[[e_1[e_1[e_2e_3]]],[e_2[e_1[e_1[e_2[e_3e_4]\ldots]]}_{\text{base case (3)}}+[e_2,\underbrace{[[e_1[e_1[e_2[e_3e_4]]]][e_1[e_1[e_2e_3]]]]}_{\text{base case (1)}}]\\
                                                       =&\ 0.
\end{align*}

\item [\textbullet] If $i_1=1$, $j_1-j_2=1$, $j_2-i_2=1$, and $i_2=1$, then $j_2=2$. This case does not occur.\vspace{0.2in}

\item [\textbullet] If $i_1=1$, $j_1=j_2\geq 4$,
\begin{align*}
&[\eiaj{i_1}{j_1},\eiaj{i_2}{j_2}]\\
=&-[[e_{j_1}\eiajthree{i_1}{j_1-2}{j_1-1},\eiaj{i_2}{j_2}]\\
=&-[\underbrace{[e_{j_1}\eiaj{i_2}{j_2}}_{\text{length decrease}},\eiajthree{i_1}{j_1-2}{j_1-1}]\\
&-[\underbrace{[\eiaj{i_2}{j_2}\eiajthree{i_1}{j_1-2}{j_1-1}]},e_{j_1}]\\
=&\ 0.
\end{align*}

\item [\textbullet] If  $i_1=1$, $j_1=j_2=3$, then $i_2=1$ or $2$. If $i_2=1$, it is trivially true. If $i_2=2$, this is the base case (2). 

\end{enumerate}\vspace{0.3in}

\item [(ii)] If $i_1\leq i_2\leq j_1<j_2$\vspace{0.2in}

\begin{enumerate}

\item[\textbullet] If $j_2-j_1\geq 2$
\begin{align*}
&[\eiaj{i_1}{j_1},\eiaj{i_2}{j_2}]\\
=&-[\eiaj{i_1}{j_1},[e_{j_2}\eiajthree{i_2}{j_2-2}{j_2-1}]\\
=&-[e_{j_2},\underbrace{[\eiaj{i_1}{j_1}\eiajthree{i_2}{j_2-2}{j_2-1}]}]\\
&-[\eiajthree{i_2}{j_2-2}{j_2-1},\underbrace{[e_{j_2}\eiaj{i_1}{j_1}}_{\text{equals 0}}]\\
=&\ 0.
\end{align*}

\item [\textbullet] If $j_2-j_1=1$, and $j_1-i_2\geq 2$
\begin{align*}
&[\eiaj{i_1}{j_1},\eiaj{i_2}{j_2}]\\
=&-[[e_{j_1}\eiajthree{i_1}{j_1-2}{j_1-1},\eiaj{i_2}{j_2}]\\
=&-[\underbrace{[e_{j_1}\eiaj{i_2}{j_2}}_{\text{length decreases}},\eiajthree{i_1}{j_1-2}{j_1-1}\\
&-[\underbrace{[\eiaj{i_2}{j_2}\eiajthree{i_1}{j_1-2}{j_1-1}]},e_{j_1}]\\
=&\ 0.
\end{align*}

\item [\textbullet] If $j_2-j_1=1$, $j_1-i_2=1$, and $i_2-i_1\geq 2$ or $j_2-j_1=1$, and $j_1=i_2$
\begin{align*}
&[\eiaj{i_1}{j_1},\eiaj{i_2}{j_2}]\\
=&\ [\eiaj{i_2-1}{j_2},\underbrace{[e_{i_2}\eiaj{i_1}{j_1}]}_{\text{length decreases}}\\
&+[e_{i_2},\underbrace{[\eiaj{i_1}{j_1}\eiaj{i_2-1}{j_2}]}]\\
=&\ 0.
\end{align*}

\item [\textbullet]  If $j_2-j_1=1$, $j_1-i_2=1$, $i_2-i_1=1$, and $i_1>1$
\begin{align*}
&[\eiaj{i_1}{j_1},\eiaj{i_2}{j_2}]\\
=&\ [\underbrace{[e_{i_1}\eiaj{i_2}{j_2}}_{\text{length decreases}},\eiaj{i_1-1}{j_1}\\
&+[\underbrace{[\eiaj{i_2}{j_2}\eiaj{i_1-1}{j_1}]},e_{i_1}]\\
=&\ 0.
\end{align*}

\item [\textbullet] If $j_2-j_1=1$, $j_1-i_2=1$, $i_2-i_1=1$, and $i_1=1$, then $j_1=3,i_2=2,j_2=4$. This is base case (3).\vspace{0.2in}

\end{enumerate}

\item [(iii)] If $i_1<j_1\leq i_2<j_2$\vspace{0.2in}

\begin{enumerate}
\item[\textbullet] If $i_1>1$
\begin{align*}
&[\eiaj{i_1}{j_1}\eiaj{i_2}{j_2}]\\
=&\ [\underbrace{[e_{i_1}\eiaj{i_2}{j_2}}_{\text{length decreases}},\eiaj{i_1-1}{j_1}]\\
&+[\underbrace{[\eiaj{i_2}{j_2}\eiaj{i_1-1}{j_1}]},e_{i_1}]\\
=&\ 0.
\end{align*}

\item[\textbullet] If $i_1=1$, $i_2\geq j_1+2$ or $i_2=j_1$
\begin{align*}
&[\eiaj{i_1}{j_1},\eiaj{i_2}{j_2}]\\
=&\ [\eiaj{i_2-1}{j_2},\underbrace{[e_{i_2}\eiaj{i_1}{j_1}]}_{\text{equals 0}}]\\
&+[e_{i_2},\underbrace{[\eiaj{i_1}{j_1}\eiaj{i_2-1}{j_2}]}]\\
=&\ 0.
\end{align*}

\item[\textbullet] If $i_1=1$, $i_2=j_1+1$
\begin{align*}
&[\eiaj{i_1}{j_1},\eiaj{i_2}{j_2}]\\
=&\ [\eiaj{i_2-1}{j_2},[e_{i_2}\eiaj{i_1}{j_1}]\\
&+[e_{i_2},\underbrace{[\eiaj{i_1}{j_1}\eiaj{i_2-1}{j_2}]}]\\
=&-\underbrace{[\eiaj{i_2-1}{j_2},\eiajthree{i_1}{j_1}{j_1+1}]}_{\text{this is case (ii)}}\\
=&\ 0.
\end{align*}

\end{enumerate}
\end{enumerate}\vspace{0.2in}

This completes the proof. Therefore $[I',I']=0$. $\square$\vspace{0.2in}

\begin{proofs}[Proof of Lemma \ref{prop1123}]
\label{proofprop1123}
\end{proofs}

Based on identities in Lemma \ref{typeccomputation}

\begin{align*}
[[e_1e_2],\one]=&\ [[e_1[e_1[e_1[e_2e_3]]]]],e_2]+[e_1[e_2[e_1[e_1[e_2e_3]\ldots]\\
=&\ \one,\\
[[e_3e_4],\one]=&\ [[e_3[e_1[e_1[e_2e_3]]]],e_4]+[e_3[e_4[e_1[e_1[e_2e_3]\ldots]\\
=&\ \one.
\end{align*}

If $k\geq 5$, then $e_k$ is commute with $e_1,e_2,e_3$, so we have $[\eij{2i+1}{j-1}{j},\one]=0$. Now if $i=1$, then

\begin{align*}
[\one,\eij{3}{j-1}{j}]=&\ [e_3,[\one\eij{4}{j-1}{j}]]\\
                                  =&\ [e_3,[e_1[e_1[\ldots[e_{j-1}e_j]\ldots]]\\
                                  =&\ 0.
\end{align*}

By the base case (2) of Proof \ref{proofIcommute}, we know $[\one,[e_1[e_2e_3]]]=0$. As for $j=4$,
\begin{align*}
&[\one,[e_1[e_2[e_3e_4]]]]\\
=&-[[e_1[e_2e_3]],[e_1[e_1[e_2[e_3e_4]]]]]]+[e_1,[[e_1[e_2e_3]][e_1[e_2[e_3e_4]]]]]\\
=&\ (\one-\one)+(\one- \one)\\
=&\ 0.
\end{align*}

Now assume $j\geq 5$. We first consider the following two brackets.
\begin{align*}
[[e_1[e_2e_3]],\ejj]=&\ [e_1,[[e_2e_3]\ejj]]\\
                                 =& -[e_1[e_3[e_2\ejj\\
                                 =&\ 0,\\
[[e_1[e_2e_3]],\eij{1}{j-1}{j}]=&-[[e_2e_3],\ejj]\\
                                                         &+[e_1,[[e_2e_3]\eij{1}{j-1}{j}]]\\
                                                      =&\ [e_3[e_2\ejj.
\end{align*}
So,
\begin{align*}
[\one,\eij{1}{j-1}{j}]=&-[[e_1[e_2e_3]],\ejj]\\
                                    &+[e_1,[[e_1[e_2e_3]]\eij{1}{2j-1}{2j}]]\\
                                 =&\ 0.
\end{align*}$\square$\vspace{0.2in}

\begin{proofs}[of Lemma \ref{centerofec}]
\label{proofcenterofec}

\end{proofs}
Clearly, $[e_1,c]=0$. Now consider $[e_{2k+1},c]$ for $k\geq 1$. There are only three summands in the right hand side of the equation which may contribute nontrivial commutators, that is, $i=k-1,k$ or $k+1$. Based on our case (3) and (4).\vspace{0.1in}

If $i=k-1$, we have the commutator equals 
\[-(n-k+1)\eia{2k-1}{2k}+\sum_{j=1}^{k-1}(-1)^{j+k-1}\eia{2j-1}{2k}.\]

If $i=k$, we have the commutator equals
\[(n-k)\eia{2k-1}{2k}+\sum_{j=1}^{k}(-1)^{j+k}\eia{2j-1}{2k}.\]

If $i=k+1$, we have the commutator equals
\[(n-k-1)\eiajthree{2k+1}{2k+2}{2k+3}-(n-k-1)\eiajthree{2k+1}{2k+2}{2k+3}.\]

The sum of these three terms is 0. Therefore, $[e_{2k+1},c]=0$ for all $k\geq 0$.\vspace{0.1in}

As for the even case: $[e_2,c]=2n[e_2,e_1]+n[e_2,[e_1[e_1e_2]]]=-2n[e_1e_2]+2n[e_1e_2]=0$. For $k>1$, similarly we have the following nontrivial case:\vspace{0.1in}

If $i=k-1$, we have the commutator equals
\[-(n-k+1)\eiaj{2k-2}{2k}+(n-k+1)\eiaj{2k-2}{2k}.\]

If $i=k$, we have the commutator equals
\[(n-k)\eiajthree{2k}{2k+1}{2k+2}-\eiajthree{2k}{2k+1}{2k+2}.\]

If $i=k+1$, we have the commutator equals
\[-(n-k-1)\eiajthree{2k}{2k+1}{2k+2}.\]

If $i\geq k+2$, we have the commutator equals

\[(-1)^{i+k-1}\eiajthree{2k}{2k+1}{2k+2}+(-1)^{i+k}\eiajthree{2k}{2k+1}{2k+2}.\]

The sum of above terms is 0. Hence, $[e_{2k},c]=0$ for all $k\geq 1$. This shows $c$ is in the center of $\fe_{C_{2n}}$.$\square$\vspace{0.2in}

\begin{proofs}[of Lemma \ref{weightvectors}]
\label{proofweightvectors}
\end{proofs}

We prove by induction on $(1)$ and $(2)$. It is clear $\phi(\bbij{1}{2})=\begin{pmatrix}E_{11}&0\\0&-E_{11}\end{pmatrix}$, and $\phi(\bbij{3}{4})=\begin{pmatrix}E_{12}+E_{22}&0\\0&-E_{21}-E_{22}\end{pmatrix}$. So $\phi([\bar{e}_1[\bar{e}_2[\bar{e}_3\bar{e}_4]]])=\phi([\bbij{1}{2}\bbij{3}{4}])=\begin{pmatrix}E_{12}& 0\\0&-E_{21}\end{pmatrix}$. This gives us the base case of (1) and (2). Now suppose $(1)$ is true for $k-1$ and $(2)$ is true for $k$. Because we know $\phi(\bbij{2k-1}{2k})=\begin{pmatrix}E_{(k-1)k}+E_{kk}&0\\0&-E_{k(k-1)}-E_{kk}\end{pmatrix}$, we have that $\phi(\bbij{2k+1}{2k+2})=\begin{pmatrix}E_{k(k+1)}+E_{(k+1)(k+1)}&0\\0&-E_{(k+1)k}-E_{(k+1)(k+1)}\end{pmatrix}$. So

\begin{align*}
&\begin{pmatrix}E_{kk}&0\\0&-E_{kk}\end{pmatrix}\\
=&\ \begin{pmatrix}E_{(k-1)k}+E_{kk}&0\\0&-E_{k(k-1)}-E_{kk}\end{pmatrix}-\begin{pmatrix}E_{(k-1)k}&0\\0&-E_{k(k-1)}\end{pmatrix}\\
=&\ \phi(\bbij{2k-1}{2k})-\phi(\beij{2k-3}{2k-1}{2k}-\ldots+(-1)^k\beij{1}{2k-1}{2k})\\
=&\ \phi(\bbij{2k-1}{2k}-\beij{2k-3}{2k-1}{2k}+\ldots+(-1)^{k-1}\beij{1}{2k-1}{2k}).
\end{align*}
Then, we have
\begin{align*}
&\begin{pmatrix}E_{k(k+1)}&0\\0&-E_{(k+1)k}\end{pmatrix}\\
=&\ \Big[\begin{pmatrix}E_{kk}&0\\0&-E_{kk}\end{pmatrix},\begin{pmatrix}E_{k(k+1)}+E_{(k+1)(k+1)}&0\\0&-E_{(k+1)k}-E_{(k+1)(k+1)}\end{pmatrix}\Big]\\
=&\ [\phi(\bbij{2k-1}{2k}+\ldots+(-1)^{k-1}\beij{1}{2k-1}{2k}),\phi(\bbij{2k+1}{2k+2})]\\
=&\ \phi([\bbij{2k-1}{2k}+\ldots+(-1)^{k-1}\beij{1}{2k-1}{2k},\bbij{2k+1}{2k+2}])\\
=&\ \phi(\beij{2k-1}{2k+1}{2k+2}+\ldots+(-1)^{k+1}\beij{1}{2k}{2k+2}).
\end{align*}

This proves (1) and (2).\vspace{0.1in}

For (3)

\begin{align*}
&\begin{pmatrix}E_{(k+1)k}&0\\0&-E_{k(k+1)}\end{pmatrix}\\
=\ &\begin{pmatrix}E_{(k+1)k}+E_{kk}&0\\0&-E_{k(k+1)}-E_{kk}\end{pmatrix}-\begin{pmatrix}E_{kk}&0\\0&-E_{kk}\end{pmatrix}\\
=\ &\phi(\bbij{2k+1}{2k})\\
   &-\phi(\bbij{2k-1}{2k}-\beij{2k-3}{2k-1}{2k}+\ldots+(-1)^{k-1}\beij{1}{2k-1}{2k})\\
=\ &\phi(\bbij{2k+1}{2k}\\
   &-\bbij{2k-1}{2k}+\beij{2k-3}{2k-1}{2k}-\ldots+(-1)^{k}\beij{1}{2k-1}{2k}).
\end{align*}

This finishes (3).\vspace{0.1in}

As for (4)

\begin{align*}
&\begin{pmatrix}0&E_{kk}\\0&0\end{pmatrix}\\
=&\ \begin{pmatrix}0&E_{(k-1)(k-1)}\\0&0\end{pmatrix}-\begin{pmatrix}0&E_{(k-1)           (k-1)}+E_{(k-1)k}+E_{k(k-1)}+E_{kk}\\0&0\end{pmatrix}\\
   &+\Big[\begin{pmatrix}E_{kk}&0\\0&-E_{kk}\end{pmatrix},\begin{pmatrix}0&E_{(k-1)(k-1)}+E_{(k-1)k}+E_{k(k-1)}+E_{kk}\\0&0\end{pmatrix}\Big],\\
=&\ \begin{pmatrix}0&E_{(k-1)(k-1)}\\0&0\end{pmatrix}-\phi(\bar{e}_{2k-1})+\Big[\begin{pmatrix}E_{kk}&0\\0&-E_{kk}\end{pmatrix},\phi(\bar{e}_{2k-1})\Big]\\
=&\ \begin{pmatrix}0&E_{(k-1)(k-1)}\\0&0\end{pmatrix}-\phi(\bar{e}_{2k-1})\\
  &+\phi(2\be_{2k-1}+[\be_{2k-3}[\be_{2k-2}\be_{2k-1}]]+\ldots+(-1)^k\beij{1}{2k-2}{2k-1})\\
=&\ \begin{pmatrix}0&E_{(k-1)(k-1)}\\0&0\end{pmatrix}\\
  &+\phi(\be_{2k-1}+[\be_{2k-3}[\be_{2k-2}\be_{2k-1}]]+\ldots+(-1)^k\beij{1}{2k-2}{2k-1}).
\end{align*}

This gives a recursive relation of $\begin{pmatrix}0&E_{kk}\\0&0\end{pmatrix}$. We can easily see that the formula in $(4)$ satisfies this  with initial condition $\begin{pmatrix}0&E_{11}\\0&0\end{pmatrix}=\phi(\be_1)$. $\square$\vspace{0.2in}

\subsection{Facts and Proofs for Type $D$ in \ref{typed}}
\label{typedproof}

\begin{proofs}[of Lemma \ref{typedcomputation}]
\label{proofsdcomputation}

\end{proofs}

We obtain this by some recursive relations. First prove the first identity with $j=i+1$.
\begin{align*}
[[e_1e_2][e_1[e_2e_3]]]=&\ [e_1[e_2[e_1[e_2e_3]]]]=-[e_1[e_1e_2]]+[e_1[e_2e_3]]=2e_1+[e_1[e_2e_3]],\\
[[e_1[e_2e_3]][e_1[e_2[e_3e_4]]]] =&-[[e_3[e_1e_2]][e_1[e_2[e_3e_4]]]]\\
=&-[[e_3[e_1[e_2[e_3e_4]]]][e_1e_2]]-[[[e_1[e_2[e_3e_4]]][e_1e_2]]e_3]\\
=&-[[e_1e_2][e_1[e_2e_3]]]+[e_3[e_4[[e_1e_2][e_1[e_2e_3]]]]]\\
=&-2e_1-[e_1[e_2e_3]]+[e_1[e_2e_3]]=-2e_1.
\end{align*} 

Now the induction step:

\begin{align*}
&[\dei{i-1}{i},\dei{i}{i+1}]\\
=&\ [[e_i\dei{i-2}{i-1}]\dei{i}{i+1}]\\
=&-[[e_i\dei{i}{i+1}]\dei{i-2}{i-1}]-[[\dei{i}{i+1}\dei{i-2}{i-1}]e_i]\\
=&-[\dei{i-2}{i-1}\dei{i-1}{i}]+[e_i[e_{i+1}[\dei{i-2}{i-1}\dei{i-1}{i}]]]\\
=&-2(-1)^{i-1}e_1=2(-1)^ie_1.
\end{align*}

If $|j-i|\ge2$, without loss of generality, assume $j-i\ge 2$, then 
\begin{align*}
&[\dei{i-1}{i}\dei{i+1}{i+2}]\\
=&-[\dei{i-1}{i}[e_{i+2}\dei{i}{i+1}]]=[e_{i+2}[\dei{i-1}{i}\dei{i}{i+1}]]\\
=&2(-1)^{i+1}[e_{i+2}e_1]=0.
\end{align*}

So by induction
\begin{align*}
&[\dei{i-1}{i}\dei{j-1}{j}]\\
=&-[\dei{i-1}{i}[e_{j}\dei{j-2}{j-1}]]=[e_{j}[\dei{i-1}{i}\dei{j-2}{j-1}]]\\
=&[e_j,0]=0.
\end{align*}

The second identity is the same as the first one when changing $e_1$ to $e_{\bar{1}}$. Next prove the third identity by induction. Some base cases:

\begin{align*}
[[\dop e_2],[e_1e_2]]=&\ [\dop e_2]-[e_1e_2],\\
[[\dop e_2][e_1[e_2e_3]]]=&-[e_2[\dop[e_1[e_2e_3]]]]+[\dop[e_2[e_1[e_2e_3]]]]\\
=&-[e_2[\dop[e_1[e_2e_3]]]]-[\dop[e_1e_2]]+[\dop[e_2e_3]],\\
[[\dop[e_2e_3]][e_1[e_2[e_3e_4]]]]=&-[[e_3[e_1e_2]][e_1[e_2[e_3e_4]]]]\\
=&-[[\dop e_2][e_1[e_2e_3]]]+[e_3[e_4[[\dop e_2][e_1[e_2e_3]]]]]]\\
=&(-1)^2([\dop[e_1e_2]]+\sum_{s=2}^3\deiajthree{s}{s}{s+1}). 
\end{align*}

By similar argument as the first identity, induction gives us 

\begin{align*}
&[\depi{i-1}{i}\dei{i}{i+1}]\\
=&-[[e_i\depi{i-2}{i-1}]\dei{i}{i+1}]\\
=&-[\depi{i-2}{i-1}\dei{i-1}{i}]+[e_i[e_{i+1}[\depi{i-2}{i-1}\dei{i-1}{i}]]]\\
=&\ (-1)^{i-1}([\dop[e_1e_2]]+\sum_{s=2}^i\deiajthree{s}{s}{s+1}).
\end{align*}

Consider

\begin{align*}
&[\depi{i-1}{i}\dei{i}{i+2}]\\
=&-[\depi{i-1}{i}[e_{i+2}\dei{i}{i+1}]]\\
=&-[e_{i+2}[\depi{i-1}{i}\dei{i}{i+1}]]\\
=&\ [e_{i+2},(-1)^{i-1}([\dop[e_1e_2]]+\sum_{s=2}^{i}\deiajthree{s}{s}{s+1})]\\
=&\ 0.
\end{align*}

Now if $j> i+2$, by induction

\begin{align*}
&[\depi{i-1}{i}\dei{j-1}{j}]\\
=&-[\depi{i-1}{i}[e_{j}\dei{j-2}{j-1}]]\\
=&-[e_{j}[\depi{i-1}{i}\dei{j-2}{j-1}]]\\
=&\ 0.
\end{align*}

If $j=i$, and $i\ge3$

\begin{align*}
&[\depi{i-1}{i}\dei{i-1}{i}]\\
=&-[\depi{i-1}{i}[e_{i}\dei{i-2}{i-1}]]\\
=&-[e_{i}[\depi{i-1}{i}\dei{i-2}{i-1}]]\\
=&\ [e_{i},(-1)^{i-1}([\dop[e_1e_2]]+\sum_{s=2}^{i-1}\deiajthree{s}{s}{s+1})]\\
=&\ \deiajthree{i-2}{i-1}{i}-\deiajthree{i-2}{i-1}{i}\\
=&\ 0.
\end{align*}

This completes the proof.
$\square$\vspace{0.5in}

\end{document}